\documentclass{article}
\usepackage{amssymb,latexsym} 
\newcommand{\1}{{{\mathchoice {\rm 1\mskip-4mu l} {\rm 1\mskip-4mu l} 
{\rm 1\mskip-4.5mu l} {\rm 1\mskip-5mu l}}}}

\newcommand{\R}{{\mathbb R}} 

\newcommand{\Q}{{\mathbb Q}} 
\newcommand{\SO}{{\rm SO}} 
 
\newcommand{\Symp}{{\rm Symp}} 
\newcommand{\id}{{\rm id}} 
\newcommand{\nnull}{{\rm null}} 
\newcommand{\sgrad}{{\rm sgrad}}

 \newcommand{\ev}{{\rm ev}} 
 \newcommand{\PD}{{\rm PD}}
  \newcommand{\oI}{{\ov{I}}}

\newcommand{\Hh}{{\mathcal H}}

\newcommand{\Z}{{\mathbb Z}} 
\newcommand{\C}{{\mathbb C}}

\newcommand{\ov}{\overline }

\newcommand{\G}{{\rm G}} 
\newcommand{\p}{{\partial}} 
\newcommand{\al}{{\alpha}} 
 
\newcommand{\be}{{\beta}} 
 
\newcommand{\Om}{{\Omega}} 
\newcommand{\om}{{\omega}}
\newcommand{\io}{{\iota}}
\newcommand{\eps}{{\varepsilon}} 
\newcommand{\de}{{\delta}} 
 
\newcommand{\ga}{{\gamma}} 
\newcommand{\Ga}{{\Gamma}} 
\newcommand{\ka}{{\kappa}} 
\newcommand{\la}{{\lambda}} 
\newcommand{\si}{{\sigma}}

\newcommand{\Ll}{{\mathcal L}}

\newcommand{\Nn}{{\mathcal N}} 
\newcommand{\Mm}{{\mathcal M}}

\newcommand{\La}{{\Lambda}}

\newcommand{\Ham}{{\rm Ham}} 
\newcommand{\area}{{\rm area\,}} 

\newcommand{\MS}{{\medskip}}

\newcommand{\NI}{{\noindent}} 
\newcommand{\proof}[1]{\noindent{\bf Proof#1:\ }} 
 
\newcommand{\QED}{\hfill$\Box$\medskip}

\newcommand{\TTilde}{\widetilde}
\newcommand{\THam}{{\TTilde{\Ham}}}
\newcommand{\trho}{{\TTilde{\rho}}}
\newcommand{\tnu}{{\TTilde{\nu}}}
\newcommand{\tphi}{{\TTilde{\phi}}}

\newcommand{\ta}{{\TTilde{a}}}

\newcommand{\CP}{{\mathbb CP}} 

\newtheorem{theorem}{Theorem}[section] 
\newtheorem{cor}[theorem]{Corollary} 
 
\newtheorem{defn}[theorem]{Definition} 
\newtheorem{thm}[theorem]{Theorem}

\newtheorem{rmk}[theorem]{Remark}
\newtheorem{lemma}[theorem]{Lemma}

\newtheorem{prop}[theorem]{Proposition}

\begin{document} 

\title{Geometric variants of the Hofer norm}
\author{Dusa McDuff\thanks{Partially
supported by NSF grant DMS0072512}\\ State University of New York
at Stony Brook \\ (dusa@math.sunysb.edu)}

\date{March 12, 2001, revised January 4, 2002 and August 30, 2006}
    
\MS\MS

\maketitle
\begin{abstract}

This note discusses some geometrically defined seminorms  on 
the group $\Ham(M, \omega)$ of Hamiltonian 
diffeomorphisms of a closed symplectic manifold $(M, \omega)$, 
giving conditions under which they are nondegenerate
and explaining their relation to the Hofer norm.  
As a consequence we show that if 
an element in $\Ham(M, \omega)$  is sufficiently close to the identity in the 
$C^{2}$-topology then it may be joined to the identity by a path whose 
Hofer length is minimal among all paths, not just among paths in the 
same homotopy class relative to endpoints.  Thus, true 
geodesics always exist for the Hofer norm. The main step in the 
proof is to show that a \lq\lq weighted"
 version of the nonsqueezing theorem holds for all fibrations 
over $S^2$ generated by sufficiently  short loops.
Further, an example is given showing that the Hofer norm  may 
differ from the sum of the one sided seminorms.
\end{abstract}
\vspace{2.8in}

\NI
\MS

\NI
keywords:  Hamiltonian symplectomorphism, 
Hamiltonian group, Hofer norm, Hofer metric, nonsqueezing 
theorem

\NI The revision of August 2006 corrects minor mistakes in the statements of
Corollary~\ref{cor:nsgeo} and Lemma~\ref{le:goods}.
\newpage

\section{Introduction}

This paper  considers some foundational 
questions about seminorms  on 
the  Hamiltonian group that were raised in 
Polterovich's lovely book~\cite{Pbk}.  
The interest of these seminorms lies in their geometric 
interpretation in terms of 
the minimal area, or equivalently curvature, of 
associated fibered spaces over the $2$-disc $D$ and
$2$-sphere $S^2$.  This allows one to  use 
geometric methods to find lower bounds for these seminorms and hence 
also for the usual Hofer norm.  
 The main question, to which we
give only a partial answer, is whether they are norms. 

Our approach gives  just one way of measuring the 
size $\rho^{+}(\phi)$ of \lq\lq one side'' 
of a Hamiltonian symplectomorphism $\phi$ but there are 
three associated two sided seminorms, the largest of which is the 
Hofer norm $\rho$.  The middle one, which we shall call $\rho_{f}$, is 
more natural from a geometric point of view, and we shall see that it 
is always nondegenerate.  However we can prove the nondegeneracy of
the smallest one, which is the sum of 
the one sided seminorms $\rho^{+}(\phi) + \rho^{+}(\phi^{-1})$, only 
in special cases, for example if $M = \CP^{n}$ or is weakly exact.

 Lalonde--McDuff~\cite{LMe} established the 
nondegeneracy of the Hofer norm for arbitrary (even open) $(M, \om)$ 
by using the nonsqueezing theorem  for the product $(M\times S^2, \om + 
dx\wedge dy)$.  In order to extend this to the more general
seminorms considered here, we need to 
understand which nontrivial symplectic  fibrations $M\to P \to S^2$
still have the nonsqueezing property.  As explained in~\S\ref{sec:sloop}
 we tackle this 
question by looking at the modified Seidel representation:
$$
\Psi: \pi_1(\Ham(M)) \to (QH_{\ev}(M))^\times,
$$
where  $(QH_{\ev}(M))^\times$ denotes the commutative group of 
multiplicative units in the 
even part of the quantum homology 
ring of $M$.

In this introduction we first define the seminorms and their geodesics,
 then give their 
geometric interpretation, and finally discuss 
short loops and nonsqueezing.   A different version of some of these results 
is being developed by Oh (see~\cite{Oh}) using the action functional
of Floer homology.  It would be very interesting to work out the relation 
between his new norm and the ones discussed here.
We will assume throughout that $(M, \om)$
is a closed connected symplectic manifold.  Further we write 
$\Ham = \Ham(M, \om)$  for the group of Hamiltonian symplectomorphisms and 
$\THam = \THam(M,\om)$ for its universal cover.
 
\subsection{The family of seminorms}\label{ss:fam}

Each  Hamiltonian
 $H_t: M\to \R, 0\le t \le 1$, defines a family of vector fields $X_t$
 (often called the {\bf symplectic gradient})
 given by\footnote
 {
 Note the sign.  Since signs are crucial when studying one sided 
 norms, we have chosen to make the signs consistent with those in 
 Polterovich's book~\cite{Pbk} even though this is different from our 
 usual conventions.
 } 
 $$
\om(X_t,\cdot) = -dH_t
 $$
 that integrates to a path
$\phi_t^{H}, 0\le t\le 1$, in $\Ham$ starting at the identity $id$.
Let 
$$
c_t = \frac 1{n!} \int_M H_t \,\om^n
$$
be the corresponding set of mean values.
We define the positive (resp. negative) length 
of $H_t$ to be:
$$
\begin{array}{lclcl}
\Ll^{+}(\phi_{t}^{H}) & = & \Ll^{+}(H_{t}) & = 
& \int_0 ^1 \left(\max_{x \in M} H_t(x) - c_t\right)\, dt,\vspace{.1in}\\
\Ll^{-}(\phi_{t}^{H}) & = & \Ll^{-}(H_{t}) & = &
\int_0 ^1 \left(c_t  - \min_{x\in M} H_t(x)\right)\, dt.\end{array}
$$
These measurements of the length of paths give rise to
seminorms on the groups $\THam$ and $\Ham$ as follows.
Recall that an element $\tphi$ of $\THam$ is
a homotopy class of paths from $id$ to 
some element $\phi \in \Ham$. 

\begin{defn}\label{def:snorm} 
We define $\trho\,^{+}(\tphi)$ (resp. $\trho\,^{-}(\tphi)$) to be the 
infimum of $\Ll^{+}(\phi_{t}^{H})$ (resp. $\Ll^{-}(\phi_{t}^{H})$)
 taken over all paths 
in the homotopy class $\tphi$.  
The corresponding seminorms $\rho^{\pm}(\phi)$ on $\Ham$
are defined by taking the infimum of $\Ll^{\pm}$ over {\bf all }
paths in $\Ham$ from $id$ to $\phi$.\footnote
{
By using suitable reparametrizations, it is not hard to see that one 
gets the same seminorms if one defines $\Ll^+(H_t)$ (resp. 
$\Ll^{-}(H_t)$) to be 
the maximum of $H_t(x) - c_t$ (resp. $c_t - H_t(x)$)
over all $t,x$: see Polterovich~\cite{Pdm}.}
\end{defn}

The seminorms 
$\nu = \trho\,^{\pm}, \rho^\pm$ have the property that for all elements 
$g,h$ of the appropriate group $\G = \THam$ or $\Ham$
$$
\nu(gh) \le \nu(g) + \nu(h),\qquad  \nu(g) = \nu(hgh^{-1}).
$$
Hence the corresponding \lq\lq distance function''
$$
d_\nu(g,h) = d_\nu(id,hg^{-1}) = \nu(hg^{-1})
$$
 is invariant under multiplication from both the left and the right.
 However, the seminorms $\nu$ under present 
consideration are not symmetric (i.e.
invariant under taking the inverse),
 though clearly
$$
\rho^{+}(\phi) \;=\; \rho^{-}(\phi^{-1}),\qquad \trho\,^{+}(\tphi) 
\;=\; 
\trho\,^{-}(\tphi^{-1}).
$$
Correspondingly $d_\nu(g,h)$ need not equal  $d_{\nu}(h,g)$.

These seminorms $\nu$ are one sided
in the sense that they only measure the size of \lq\lq half" of $H_t$,
either the part $H_t - c_t$ above the mean, or the part below.
We now discuss several different  ways to obtain
 two sided norms.  By definition these are  symmetric,
i.e. invariant under taking the inverse, so that the corresponding 
pseudometrics will be biinvariant.

The usual Hofer norm $\rho$ on $\Ham$ is given by 
minimizing\footnote
{
By slight abuse of language we often use the word \lq\lq minimizing'' 
in situations when the minimum (or more properly infimum) is not attained.}
the sum 
$$
\Ll(H_t): = \Ll^+(H_t) + \Ll^-(H_t);
$$
we write
$\trho$ for the associated seminorm on $\THam$.  Thus $\trho(\tphi)$
is the infimum of $\Ll^+(H_t) + \Ll^-(H_t)$ over all paths in the homotopy 
class  $\tphi$, while $\rho(\phi)$ is its infimum over all paths from $id$ 
to $\phi$.  Clearly,
$$
\trho(\tphi) = \trho(\tphi^{-1}),\qquad \rho(\phi) = \rho(\phi^{-1}).
$$
Although $\rho$ is known to be a norm for all $M$ (cf~\cite{LMe}), 
it is not clear 
whether $\trho$ is always a norm.  Its null set
$$
\nnull(\trho) = \{\tphi: \trho(\tphi) = 0 \}
$$
is a normal subgroup, and so must lie in the kernel $\pi_1(\Ham)$ of the
covering map $\THam \to \Ham$.\footnote
{
Recall that $\Ham(M)$ is a simple group when  $M$ is closed.
Thus a symmetric, conjugation 
invariant seminorm on $\Ham(M)$ must be either  
nondegenerate or everywhere zero.}
It is conceivable that 
$\nnull(\trho)$ is nonempty for some $M$.  For example, 
there might be 
nontrivial elements of 
$\pi_1(\Ham(M, \om))$ that are supported in a Darboux chart.
Such elements lie in $\nnull(\trho)$ since the size of their 
generating Hamiltonian can be made arbitrarily small by a suitable 
conformally symplectic conjugation.  They do not exist when 
when $\dim(M) \le 4$  since the group of compactly supported 
symplectomorphisms of Euclidean space is 
contractible in these dimensions.  
However, they might exist in higher dimensions.  In any case, $\trho$ 
descends to a norm on the quotient group $\THam/\nnull(\trho)$.

We will also consider the two sided seminorm on $\THam$
given by the
sum  $\trho\,^{+} + \trho\,^-$.  
It is easy to see that
$$
\trho(\tphi)\;\; \ge\;\; \trho\,^{+}(\tphi) + \trho\,^-(\tphi).
$$
In principle, there could be strict inequality here
 since the one sided seminorms 
$\trho\,^{+}(\tphi)$ and $\trho\,^{-}(\tphi)$
could well be realized by different minimizing 
sequences of Hamiltonians.  However we have so far found no example 
to illustrate this possibility.\footnote
{In~\cite{En} Entov considers the seminorm $\max(\trho\,^{+},\trho\,^-)$.
This is clearly equivalent to $\trho\,^{+} + \trho\,^-$, and it is easy to 
check that it has the same 
geodesics.}
Sometimes we will write $\trho_f$ instead of $\trho^+ + \trho^-$.

There are two possible ways to obtain a related seminorm on $\Ham$;
one may either  consider 
the  sum $\rho^{+} + \rho^-$ or consider
 the seminorm $\rho_f$ on $\Ham$ induced by $\trho\,^{+} + \trho\,^-$, viz:
$$
\rho_f = \inf \{\trho\,^{+}(\tphi) + \trho\,^-(\tphi):
\tphi \mbox{ lifts } \phi\}.
$$
Both $\rho^{+} + \rho^-$ and $\rho_f$  are symmetric.  Further 
$$
\rho^{+}(\phi)  + \rho^-(\phi) \;\;\le \;\;\rho_f(\phi) \;\;\le\;\; \rho(\phi).
$$

Polterovich asked in~\cite{Pbk} whether the seminorm
$\rho^+  + \rho^-$  must always be a norm, and indeed whether 
it always equals $\rho$. 
The  proof given in~\cite{LM} that $\rho$ is a norm adapts almost immediately to
show that $\rho_f$ is always a norm.  However, the question for
$\rho^+  + \rho^-$ is much more difficult and revolves around 
properties of $\pi_1(\Ham)$ and the quantum homology of $M$
that are discussed in more detail in \S\ref{ss:sloop}. 
We have only succeeded 
in showing that it is a 
norm in special cases, some of which are 
mentioned in the next result. Recall that $(M, \om)$ is said to be 
{\bf weakly exact} if $\om$ vanishes on the image $H^{S}_{2}(M, \Z)$ of $\pi_2(M)$
in $H_{2}(M, \Z)$. 

\begin{thm}\label{thm:norm} \begin{itemize}
\item[(i)] $\rho_f$ is a norm for all closed $M$.  

\item[(ii)] $\rho^+  + \rho^-$ is a norm if $M$ is weakly exact or 
$M = \CP^{n}$.

\item[(iii)] There is a symplectic form $\om$ on the one point blow up $M_*$ of 
$\CP^2$ such that $\rho^+  + \rho^-$  is a norm
on $\Ham(M_*, \om)$ that is distinct from $\rho_f$ and hence also from 
$\rho$.\end{itemize}
\end{thm}

We also show that $\rho^+  + \rho^-$ is a norm whenever there are no 
asymmetric short loops in the sense of Definition~\ref{def:sloop} below. 
Other conditions under which $\rho^+  + \rho^-$ is a norm are 
given later.  For example, it would suffice that
the nonsqueezing theorem hold for  all  $[\la]\in 
\pi_{1}(\Ham(M, \om))$: see Corollary~\ref{cor:ndeg} and Lemma~\ref{le:nsgeo}.
  The example in (iii) is discussed in more detail in 
\S\ref{ss:exam}.

The question of whether  one sided seminorms  such as $\rho^+$ 
are nondegenerate seems intractable.  Since the corresponding   
 null set
is a conjugation invariant semigroup rather than a normal subgroup,
the simplicity of the group $\Ham(M, \om)$ is now of no help: 
to prove nondegeneracy 
one must find a lower bound for $\rho^+(\phi)$ for
 {\it every} element $\phi\in \Ham(M,\om)$, a task beyond the reach 
 of the  techniques used here.

 \begin{rmk}[Noncompact manifolds]\rm  The above definitions are 
generalizations of notions first introduced 
for Euclidean space $(\R^{2n}, \om_0)$.  More generally, given any
 noncompact manifold $(M, \om)$ (without boundary) let
 $\Ham^c(M, \om)$ be the group of compactly supported Hamiltonian 
 symplectomorphisms and consider Hamiltonians $H$ with compact support.  
 Then set
 $$
 \Ll^+(H_t): = \int \max(H_t)\,dt,\qquad  \Ll^-(H_t): = -\int 
 \min(H_t)\,dt,
 $$
 and use these length measurements to define seminorms just as before.
 Clearly, in this case 
 there are  elements $\phi \ne id$ such that $\rho^+(\phi) = 0$, 
for example the time $1$ map of a nonzero function $H_t$ with $H_t \le 0$.
One can recover a situation more like the closed case by restricting 
attention to the subgroup $\Ham_0^c(M, \om)$ generated by compactly 
supported Hamiltonians of zero mean.  This subgroup, the kernel of 
the Calabi homomorphism, is simple just as in the closed case, and even 
in the case $M = \R^{2n}$ it is not yet known whether the seminorms 
$\rho^{\pm}$ are nondegenerate on it.

On the other hand  a few of the questions considered here are more 
tractable when $M = \R^{2n}$. For example, the sum $\rho^+ + \rho^-$ 
is always nondegenerate.  This was first proved by Viterbo~\cite{Vit}
who used generating functions to construct a section $\phi\mapsto 
c(\phi)$ of the action spectrum bundle with the property that
$$
0 < c(\phi) \le \rho^+(\phi) + \rho^-(\phi),\quad\mbox{when 
}\;\;\phi\ne id.
$$
 It also follows 
from our arguments since the fact that the
elements of $\pi_1(\Ham^c(\R^{2n}, \om_0))$
have representatives with arbitrarily short Hofer length implies that
$\rho^+ + \rho^- = \rho_f$: see~\S1.4.   However, it is still unknown 
whether the two norms $\rho^+ + \rho^-$ and $\rho$ must always agree, 
even when $M = \R^2$ or $\R^4$.  Note also that when $M$ is Euclidean 
space there are several
other possible one sided seminorms arising from various selectors and 
it is not clear that $\rho^{\pm}$ are the 
most interesting ones to consider.  For further discussion see
Polterovich~\cite{Pbk} and  Schwarz~\cite{Sch}.
\end{rmk}

\subsection{Results on Geodesics}\label{ss:geod}

We now  describe our results on geodesics.  It is possible to define 
geodesics as critical points of a suitable length functional: 
see~\cite{Pbk,LM}.  However 
the lack of smoothness of this functional causes some problems.  We will 
take a different approach that nevertheless gives rise to the same geodesics.
We write $\G$ to denote one of the groups $\THam$, $\Ham$.

Observe that the seminorms under consideration are of two kinds.
If $\nu$ is one of $\trho^\pm, \trho, \rho^\pm$ and $\rho$ then $\nu(g)$
is the infimum of
an appropriate length functional $\Ll_\nu$ over all paths in $\G$ 
from $id$ to $g$, while the value of the other seminorms  $\trho_f, \rho_f$ and 
$\rho^+ + \rho^-$ at an element $g$ depend on
 minimizing the functionals 
$\Ll^\pm$ on two possibly different paths.  If $\nu$ is one of the three
latter seminorms we set $\Ll_\nu = \Ll$. 

\begin{defn}  Let $(\nu, \Ll_\nu)$ be one of the pairs defined above.
A path  $\{g_{t}\}_{t\in 
[a,b]}$ in $\G$  is said to be {\bf $\nu$-minimizing} if it achieves the minimum of 
$\Ll_\nu$, i.e. if
$$
\nu(g_b g_a^{-1}) = d_\nu(g_a, g_b) = \Ll_\nu(\{g_t\}).
$$
When defining geodesics, it is convenient to restrict attention to 
paths $\{g_{t}\}$ that are {\bf regular} in the sense that their
generating vector field $\dot{g}_{t}$ is never zero.
We then say that a regular
path $\{g_{t}\}_{t\in 
[a,b]}$ is a {\bf $\nu$-geodesic} from  $g_{a}$ to $g_{b}$ if
there is $\eps > 0$ such that for all $t_{0}\in [a,b]$ each path 
$$
[a,b]\cap [t_{0}, t_{0} +\eps] \to \G:   \quad
t\mapsto g_{t}g_{t_{0}}^{-1}
$$
is $\nu$-minimizing.  
\end{defn}

It follows immediately that any $\rho_{f}$ or $(\rho^+ +\rho^-)$-geodesic 
 is also a $\rho$-geodesic.
Observe also this definition  is somewhat 
stronger than is usual  in this context (cf. \cite{LM,Pbk})
since it requires that short 
pieces of $\nu$-geodesics minimize $\Ll_\nu$ among {\it all} paths
 with the given endpoints rather than 
simply among the homotopic paths.   Though natural, 
it was not used previously for the group $\G = \Ham(M)$
 because it was not known
 whether there always are paths 
satisfying this stronger condition.  This is
 the problem of \lq\lq short loops'' that is discussed in
 \S\ref{ss:sloop}.  The main new step
in the  proof of Theorem~\ref{thm:geod} below is to 
get around this difficulty.

Before stating it, we recall some ideas from Bialy--Polterovich~\cite{BP}
and Lalonde--McDuff~\cite{LM}.  A path 
$\ga = \{\phi_{t}\}_{t\in [a,b]}$ in $\Ham$ that is
generated by the Hamiltonian $H_{t}$ 
is said to have a {\bf fixed maximum (minimum)} 
if there is a point $x_{0}\in M$ such that
$$
H_{t}(x)\;\;\le (\ge)\;\; H_{t}(x_{0}),\quad\mbox{for all}\;\; x\in M, \;t\in [a,b].
$$
It is said to have a fixed maximum (minimum) {\bf at each moment} if there 
is $\eps > 0$ such that each subpath 
$\{\phi_{t}\phi_{t_{0}}^{-1}\}_{t\in [a,b]\cap [t_{0}, t_{0} +\eps]}$ has a fixed 
maximum (minimum).  
It was shown in~\cite{BP} that
a path $\ga$ in $\Ham^{c}(\R^{2n})$ is a $\rho$-geodesic
if and only 
if it has both a fixed maximum and a fixed minimum at each moment.
The result proved in~\cite{LM} for general $M$ can be stated in our 
current language as follows:\footnote{
 The papers Lalonde--McDuff~\cite{LM} 
were written before 
it was understood how to define Gromov--Witten invariants on arbitrary 
symplectic manifolds.   Therefore,
many of the results in part II have unnecessary restrictions. 
In particular, 
in Theorems 1.3 (i) and 1.4 and in Propositions 1.14 and 1.19 (i)
one can remove the hypothesis that
$M$ has dimension $\le 4$ or is semi-monotone since the proofs are based 
on the fact that
 quasicylinders $Q = (M\times D^{2}, 
\Om)$ have the nonsqueezing property.  The result we now quote 
from~\cite{LM} incorporates these improvements. All the details of 
 its proof (besides the construction of 
general Gromov--Witten invariants) occur in~\cite{LM}: see~\S\ref{sec:geod} 
below.  Oh gives a new proof in~\cite{Oh}.}    
a path $\ga$ in $\THam(M, \om)$ is a 
$\trho$-geodesic if and only 
if it has both a fixed maximum and a fixed minimum at each moment.
This generalizes to the one sided seminorms on 
$\THam$ as follows.

\begin{prop}\label{prop:geod}
 Let $\ga = \{\phi_{t}\}_{t\in [a,b]}$ be a 
path in $\Ham(M)$. 
\begin{itemize}
\item[(i)]
A lift of $\ga$ to $\THam$ 
is a geodesic with respect to $\trho\,^+$ (resp.
$\trho\,^{-}$) if and 
only if $\ga$ has a fixed minimum (resp. maximum) at each moment;
\item[(ii)]  A lift of $\ga$ to $\THam$ 
is a geodesic with respect to
 $\trho\,^{+}+\trho\,^{-}$ if and only 
    if $\ga$ has both a fixed maximum and a fixed minimum at each moment.
\end{itemize}
\end{prop}

It is much harder to get results on $\Ham(M)$. 
Here is our main result.

\begin{thm}\label{thm:geod}  Let $\ga = \{\phi_{t}\}_{t\in [a,b]}$ be a 
path in $\Ham(M)$. 
\begin{itemize}
\item[(i)] If $\nu $ is one of the norms $\rho$ or $\rho_f$
on $\Ham(M)$, then $\ga$ is a $\nu$-geodesic  
if and only 
    if it has both a fixed maximum and a fixed minimum at each moment.
\item[(ii)]  If $M = \CP^{n}$ or is weakly exact,
 the same result holds 
with $\nu = \rho^+ + \rho^-$.
\end{itemize}
\end{thm}

\begin{cor}\label{cor:ndeg}  $\rho_f$ is nondegenerate for all $M$, 
while $\rho^+ + \rho^-$
is nondegenerate  when $M$ satisfies the conditions in (ii) above.
\end{cor}
\proof{}  Because $\Ham(M)$ is a simple group, a symmetric conjugation invariant 
seminorm  $\nu$ on $\Ham(M)$ is either identically zero or is 
nondegenerate.  But in the 
former case every path in $\Ham$ would be a geodesic.  Hence,
 the existence of nongeodesic 
paths implies that $\nu$ is nondegenerate.\QED

The proof of the above theorem finds lower bounds for 
$\nu(\phi)$ when $\phi$ is close to the identity
by using the geometric characterization 
of the seminorms given in \S\ref{ss:geomint} and suitable extensions of 
the nonsqueezing theorem.  
Further, we generalize~\cite{BP,LM} by proving the following 
local flatness result for a neighborhood of $id$ in $\Ham(M, \om)$.  
Recall that $\Ham(M, \om)$ is the kernel of a surjective homomorphism
$$
{\rm Flux}:\quad \Symp_{0}(M, \om)\to H^{1}(M, \R)/\Ga_{\om}
$$
where the flux group $\Ga_{\om}$ is finitely generated but is not 
known to be discrete in all cases. 
  Hence the most we can say in general
  is that $\Ham(M, \om)$ sits inside the identity 
component $\Symp_{0}$ as the leaf of a foliation.  Therefore, we do 
not use the topology on $\Ham$ induced from $\Symp_{0}$ but instead 
use the topology on $\Ham$ induced from the $C^{2}$-topology on the Lie algebra 
of Hamiltonian functions with zero mean.   Thus a neighborhood of the 
identity  consists of all time $1$-maps of Hamiltonian flows 
generated by Hamiltonians $H_{t}$ that are sufficiently small in the 
$C^{2}$ topology.

\begin{prop}\label{prop:nbhd}  There is a path connected 
neighborhood $\Nn\subset \Ham(M)$ of the identity  in the $C^{2}$-topology
such that any element $\phi\in \Nn$ can be joined to the identity by 
a path that minimizes both $\rho$ and $ \rho_f$.  In 
particular, these two norms agree on $\Nn$.  Moreover, $(\Nn, \rho)$ 
is isometric to a neighborhood of $ \{0\}$ in a normed vector space.
\end{prop}

 One can also look for longer
length minimizing paths starting from $id$.  Note that one cannot reach an arbitrary
$\phi\in \Ham(M)$ by such a path: an example is given  in~\cite{LM}~II of an 
element in $\Ham(S^{2})$ that cannot be reached by any 
$\rho$-minimizing path.  Moreover, since $(\rho^+ + \rho^-)$-geodesics 
are the same as $\rho$-geodesics (when they exist), any path that minimizes 
$\rho^+ + \rho^-$ also minimizes $\rho$.  Hence  elements $\tau$ 
with $0\ne \rho^+(\tau) + \rho^-(\tau) < \rho(\tau)$ as in 
Theorem~\ref{thm:norm}(iii) cannot be reached by   $(\rho^+ + 
\rho^-)$-minimizing paths.\footnote
{
In fact, such an element $\tau$ cannot be reached even by a 
$(\rho^+ + \rho^-)$-minimizing
sequence $\ga_i$ of paths, i.e. such that
$\rho^+(\tau) + \rho^-(\tau)$ is the limit of $\Ll^+(\ga_i) + 
\Ll^-(\ga_i)$.  Similarly, if there were an element $\phi$ such that 
$\rho_f(\phi) \ne \rho(\phi)$ then it could not be reached by a 
$\rho_f$-minimizing sequence of paths.}

The following result is a mild extension of 
work by Entov~\cite{En} and McDuff--Slimowitz~\cite{MSlim}. 
(See also Oh~\cite{Oh} where the result is generalized to some time 
dependent paths.)
We will say that a time independent Hamiltonian $H$ is {\bf slow} if
neither its flow  nor the linearized flows at its critical points
have nonconstant contractible periodic orbits of period $ <1$.

\begin{prop}\label{prop:good}  Let $H$ be a slow Hamiltonian.  Then
the path $\phi_{t}^{H}, 0\le t\le 1$ minimizes both $\trho\,^{-}$
and $\trho\,^{+}$
and hence also minimizes $\trho$ on $\THam(M)$.  
If in addition $(M,\om)$ is weakly 
exact, this path minimizes all the norms
$\rho^+ + \rho^-$, $\rho_f$ and $\rho$ on $\Ham(M)$.
\end{prop}

\subsection{Geometric interpretations of the seminorms}\label{ss:geomint}

 Consider a smooth fibration $\pi: P  \to 
B$ with fiber $M$, where $B$ is either  $S^2$ or the $2$-disc $D$.
Here we consider $S^2$ to be the union $D_+ \cup D_-$ of two copies of 
$D$, with the same orientation  as $D_+$.  We denote the equator $ D_+\cap 
D_-$ by $\p$, oriented as the boundary of $D_+$, and choose some point
$*$ on $\p$ as the base point of $S^2$.  Similarly, $B = D$ is provided 
with a basepoint $*$ lying on $\p = \p D.$  In both cases, we assume that 
the fiber over $*$ has a chosen identification with $M$.

In this paper we will be considering triples $(P,\pi,\Om)$ where 
$\pi:P\to B$ is a fibration as above and $\Om$ is a
normalized $\pi$-compatible
 symplectic form  on  $P$.  $\pi$-compatibility
 means that the restriction $\om_b$
of $\Om$ to the fiber $M_b = \pi^{-1} (b)$ is nondegenerate for each $b\in B$,
and the normalization condition is that $\om_* = \om$. 
For short we will write $(P, \Om)$ instead of the triple $(P,\pi,\Om)$, and 
will use the words \lq\lq $\om$-compatible" instead of 
\lq\lq normalized $\pi$-compatible."  

In \S\ref{ss:area} we will describe in more detail 
the geometric structure that is induced by $\Om$ on the fibration
$\pi:P\to B$.   The most important point is that $\Om$ defines a
connection on $\pi$ whose horizontal distribution is $\Om$-orthogonal 
to the fibers.   
If $\al$ is any path in $B$ then $\pi^{-1}(\al)$ is 
a hypersurface in $P$ whose characteristic foliation consists of 
the horizontal lifts of $\al$,
and it is not hard to check that
the resulting holonomy is Hamiltonian round every
contractible loop.   Because $B$ is simply connected, it follows 
that the structural group of $\pi$ can be reduced to 
 $\Ham(M, \om)$.  Further, $\pi$ can be symplectically 
trivialized over each disc $D$
 by parallel translation along 
a suitable set of  rays.  This means that there is a fiber preserving 
mapping  
$$
\Phi:  \pi^{-1}(D) \to M\times D, \quad \Phi|_{M_*} = id_M
$$
such that the pushforward $\Phi_*\Om$ 
 restricts to the same form $\om$ on each fiber $M\times pt$.
 The fibration $(P, \Om)\to S^2$ is said to be 
{\bf symplectically trivial}\footnote
{
The words \lq\lq symplectically trivial"  mean 
\lq\lq trivial as a symplectic fibration."  In the present context,
this implies  triviality as a Hamiltonian fibration.} 
if such a map exists from $P$ to $ M\times S^2$.

\begin{defn}
The {\bf monodromy}  $\phi = \phi(P)\in \Ham(M)$  
of a fibration $(P, \Om)\to B$
is defined to be the monodromy of the connection determined by 
$\Om$ around the based loop $(\p, *)$.
 Using the trivialization of
$P$ over $\p$ provided by $B$ itself if $B = D$
or by $D_+$ if $B = S^2$, one gets  a well defined lift 
$\tphi$ of $\phi$ to $\THam$.  Sometimes we will write $P_\tphi$ 
(resp.  $P_\phi$) for a 
fibration $(P, \Om) \to B$ with monodromy $\tphi$ (resp. $\phi$).
\end{defn}

The next definition describes various different  area measurements.

\begin{defn}\label{def:area}
The {\bf area} of a fibration $(P, \Om) \to B$ 
 is defined to be:
$$
{\rm area}\,(P, \Om) = \frac {{\rm vol\,}(P,\Om)}{{\rm vol\,}(M, \om)}
= \frac{\int_{P}\Om^{n+1}}{(n+1)\int_{M}\om^{n}}. 
$$
Further: \begin{itemize}
\item[(i)]
 $\ta^+(\tphi)$ (resp. $a^+(\phi)$)
is the infimum of $\area (P, \Om)$ taken over 
all $\om$-compatible symplectic forms $\Om$ on the fibration $P\to D$
 with monodromy $\tphi$ (resp. $\phi$).
\item[(ii)]  $\ta(\tphi)$
(resp. $a_f(\phi)$) is 
the infimum of $\area(P, \Om)$ taken over all symplectically trivial 
fibrations $(P, \Om)\to S^2$ with monodromy $\tphi$ (resp. $\phi$).
\item[(iii)]
$a(\phi)$ is the infimum of $\area(P, \Om)$ taken over {\it all} 
fibrations $(P, \Om)\to S^2$ with monodromy $\phi$.
\item[(iv)]  $\ta^-(\tphi) = \ta^+(\tphi^{-1})$ and 
$a^-(\phi) = a^+(\phi^{-1})$
\end{itemize}
\end{defn}

It is easy to see that  $a^+(\phi)$  (resp. $a_f(\phi)$)
is the infimum of $\ta^+(\tphi)$ (resp. $\ta(\tphi)$)
over all lifts $\tphi$ of $\phi$ 
to $\THam$.   The following lemma amplifies
Polterovich's results in~\cite{Peg}. 

\begin{prop}\label{prop:polt}\begin{itemize}
\item[(i)]  $\trho\,^{+}(\tphi) =
\ta^+(\tphi)$;
\item[(ii)]  $\rho^{+}(\phi)  + \rho^-(\phi) = a(\phi)$;
\item[(iii)]  $\rho_f(\phi) = a_f(\phi)$.
\end{itemize}  
\end{prop}

This is proved in  \S\ref{ss:area}.  We have here interpreted
 our seminorms in terms of area since this is what 
our methods estimate.  However, as is clear from the proof of
the above Proposition, these area measurements are 
equivalent to suitable measurements of curvature: see Polterovich's remarks 
about K-area in~\cite{Pk,Peg} and also Entov~\cite{En}.  

 The Hofer norm also has a 
geometric interpretation:  $\rho(\phi)$ is the infimum 
of area$(P_\phi, \Om)$ taken over
all  fibrations $(P_\phi, \Om)$
for which there is a symplectomorphism 
$$
\Phi:
(P_\phi, \Om)\to(M\times S^2, \om + \al) 
$$
 that is the identity on the fiber over the base point and takes
the hypersurface $\pi^{-1}(\de)$ lying over the equator
to a specially situated
hypersurface $\Ga_H$ in the product.  
(The precise condition on $\Ga_H$ is described at the end of 
\S\ref{ss:area}.  
Note that $\Phi$ need  not 
preserve the given fibered structure on $P_\phi$.)
Even if we ignore the condition on $\Ga_H$, we may be minimizing over a 
smaller set than in the definition of $a_{f}$, since it is not clear whether every 
$\om$-compatible form on the trivial fibration $M\times S^2 \to S^2$ is 
symplectomorphic to a product form.
From a geometric point of view the minimizing sets used to define
 $a(\phi)$ and $a_f(\phi)$ are much more 
natural than that used to define the Hofer norm. Thus 
 $\rho^+  + \rho^-$ and $\rho_f$
are the seminorms with the most geometric meaning.  We consider 
$\rho_f$ to be the geometric analog of the Hofer norm
and hence have called it $\rho_f$ where $f$ denotes \lq\lq fibered." 
Correspondingly we will sometimes write
$\trho_f$ instead of $\trho\,^+ + \trho\,^-$ for its lift to $\THam$.
\MS

Finally, let us consider the case when $\tphi$ is a lift of $id$; in other 
words, $\tphi$ is a homotopy class $[\la]$ of loops in $\Ham$.
Above we have measured $\trho^+([\la])$ by means of the area of a 
fibration over $D$ with boundary monodromy $[\la]$.
On the other hand, the natural geometric object associated to a loop
is a fibration over $S^2$ constructed by using 
the loop $\la$ as a clutching function:  
$$
P_\la = \left(D_+^2\times M\right) \cup_\la \left(D_-^2\times 
M\right),
$$
where\footnote
{
Note the direction of this attaching map:  
this is a different convention from~\cite{LMP,Mcq} though the same 
as~\cite{Pbk}.}  
$$
\la: (2\pi t, x)_-\mapsto (2\pi t, \la_t(x))_+.
$$
By fixing an
identification of the fiber of $P_\la$ at the basepoint $*\in \p D_+$
with $M$, we can normalize the loop
$\la = \{\la_t\}$ by requiring
that $\la_* = id$.
It is not hard to see that the symplectic form on the fibers
has a closed extension to $P_\la$ precisely when
$\la$ is homotopic to a loop in  $\Ham(M,\om)$. 
Thus there is a bijective correspondence between classes $[\la] \in
\pi_1(\Ham(M,\om))$ and  Hamiltonian  fibrations $P_\la$ with one 
fiber identified with $M$.  Hence, the most obvious geometric way to measure
the size of $[\la]$ is  to take the minimum area of $\om$-compatible 
forms $\Om$ on this fibration $P_\la \to S^2$.  
The next lemma says that 
this gives nothing new.

\begin{lemma}\label{le:polt2}  For each $[\la]\in \pi_1(\Ham)$, 
$\trho^+([\la])$ is the infimum  of {\rm area}$(P_{\la}, \Om)$, where 
$\Om$ ranges over all $\om$-compatible forms on the
fibration $P_{\la} \to S^2$ constructed above.
\end{lemma}

Observe here that even though $P_{\la} \to S^2$ is a fibration over the 
closed manifold $S^2$ the corresponding area measure is the
 {\it one sided} seminorm $\trho^+([\la])$.  One way to understand this is 
 to consider $P_{\la}$ as made by gluing together two fibrations, the first
$(P_+, \Om)\to D_+$
 with monodromy $\la$ and area $\ge \trho^+([\la])$
 and the second $(P_-, \Om)\to D_-$ with trivial monodromy 
 and hence arbitrarily small 
 area.

\subsection{Short loops and the Nonsqueezing theorem}\label{ss:sloop}

Though one can find lower bounds for the 
seminorms $\trho, \trho\,^{\pm}$ on $\THam$ by the methods of~\cite{LM,En,MSlim} 
it is usually more 
difficult to find lower bounds  on $\Ham$.
There is one lower bound for $\rho(\phi)$ that is independent of the path
from $id$ to $\phi$, namely the energy--capacity inequality:
\begin{quote}{\it Let $B = \phi(B^{2n}(r))$ be a symplectically
    embedded ball in $(M, \om)$ of  radius $r$.
    If $\phi(B)\cap B = \emptyset$, then}  $\rho(\phi) \ge \pi r^{2}/2$.
\end{quote}
This was proved for $\R^{2n}$ by 
Hofer~\cite{Ho} (without the constant $1/2$),
 and for any symplectic manifold in~\cite{LMe}.  
It follows immediately that $\rho$ is a norm.
However, 
this inequality does not hold for $\rho\,^{\pm}$: one 
can find counterexamples by adapting the construction in
Eliashberg--Polterovich~\cite{EP} that shows
 the degeneracy of the $L^{p}$-metrics 
for $p < \infty$.

The paper~\cite{LM}~II
proposes another way to get lower bounds for $\trho(\tphi)$ in the case 
when $\tphi$ is 
$C^{2}$-close to $id$.  The idea is to 
define the \lq\lq graph'' $\Ga_{H}$ of a Hamiltonian isotopy and to 
embed symplectic balls $B^{2n+2}(\eps)$ of radius $\eps$ both \lq\lq 
under'' and \lq\lq over'' this graph.  If $\phi_{t}$ is sufficiently 
$C^{2}$-close to the $id$ and has fixed extrema, one can 
construct such embeddings with
$$
\pi\eps^{2} = \Ll(\phi_{t}).
$$
A simple argument using the nonsqueezing theorem then implies that
the path $\{\phi_{t}\}_{t\in [0,1]}$ in $\THam(M)$  minimizes $\trho$.
(This argument is explained in more detail in \S\ref{sec:geod} below.)

As remarked in~\cite{LM}, in order to go from here to 
estimates on the group $\Ham$ we need to understand the {\bf short 
loops}.  More precisely, given a seminorm $\tnu$ on $\THam(M)$ 
define 
$$
\ell_{\tnu}: \pi_{1}(\Ham(M, \om)) \to [0, \infty)
$$
to be the restriction of $\tnu$ to the subgroup 
$\pi_{1}(\Ham)\subset \THam$.
To simplify the notation, we will denote its value on the homotopy class 
$[\la]$ of a loop $\la$ in $\Ham$ by $\ell_\tnu(\la)$.
Further, we define  $r_\tnu(M)$ to
be  the minimum of 
the {\it positive} values of $\ell_{\tnu}$.  If $\ell_\tnu(\la) = 0$
for all $[\la]$ then we set $r_\tnu(M) = \infty$.
For short, we will write
$$
\begin{array}{lllll}
\ell^\pm, r^\pm  & \mbox{for} &  \ell_\tnu, r_\tnu &\mbox{when} & \tnu = 
\trho^\pm.
\end{array}
$$
Further, we define $r_a(M)$ to be the supremum of $\de \ge 0$ such that
$$
\left(\ell^-(\la) < \de \mbox{ or } \ell^+(\la) < \de \right)\;\;
\Longrightarrow \;\; \left(
 \ell^-(\la) = 0 = \ell^+(\la).\right)
$$
Thus $r_a(M) = 0$ if there is a sequence of loops
$\la_i$ in $\Ham(M)$ such that $\ell^-(\la_i) \to 0$ 
while $\inf_i \ell^+(\la_i) > 0$.  In particular, $r_a(M) = 0$ if
 there is a loop $\la$ such that
$\ell^-(\la) = 0$ while $\ell^+(\la) > 0$.

\begin{defn}\label{def:sloop}
The manifold $M$ is said to have 
{\bf  $\tnu$-short loops} if $0$ is not an isolated point in the image
of $\ell_{\tnu}$, or equivalently if $r_{\tnu}(M) = 0$.
It is said to have {\bf asymmetric short loops} if 
$r_a(M) = 0$.  
\end{defn}

The next lemma shows that $\tnu$-geodesics in $\THam$ 
descend to $\nu$-geodesics in $\Ham$ when 
there are no (asymmetric) short loops.  

\begin{lemma}\label{le:asloop} \begin{itemize}\item[(i)] Let $(\tnu, \nu)$ be one of the pairs
$(\trho,\rho), (\trho_f,\rho_f)$.
Suppose that the path $\ga$ has $\tnu$-length $\le r_{\tnu}(M)/2$ and 
minimizes $\tnu$ in $\THam(M)$. Then it is 
$\nu$-minimizing in $\Ham(M)$.
\item[(ii)] If $(\tnu, \nu) = (\trho_f, \rho^+ + \rho^-)$
the same statement holds with $r_{\tnu}(M)$ replaced by $r_a(M)$.
\end{itemize}
\end{lemma}
\proof{} (i)   We know that $\ga$ minimizes $\Ll$ among all homotopic 
paths and need to see that it minimizes $\Ll$ among all paths with 
the same endpoints.  If not, there is another shorter path $\ga'$ 
with the same endpoints.  Then $\la = (-\ga') * \ga$ is a loop
with length
$$
\Ll(\la) \;\;\le\;\; \Ll(\ga) + \Ll(-\ga')\;\; <\;\; r_{\tnu}(M).
$$
Hence it has zero length.  This means that we can compose $\ga'$ with 
an arbitrarily $\nu$-short loop homotopic to $\la$ to obtain a path 
that is homotopic to $\ga$ but shorter than it.  This contradiction proves the 
lemma.

In case (ii) essentially the same argument works.  Now we may only assume 
 that
$\ga'$ reduces one of the seminorms $\trho\,^\pm$, say $\trho\,^-$.
Then 
$$
\trho\,^+(\la)\;\; = \;\;\trho\,^-(-\la) \;\;=\;\; 
\trho\,^-((-\ga)* \ga') 
\;\;\le \;\; \Ll\,^+(\ga) + 
\Ll\,^-(\ga') \;\;<\;\; r_a(M).
$$ 
Hence both $\trho\,^+(\la)$ and $\trho\,^-(\la)$ are $0$ and the argument 
proceeds as before.
\QED

In the general case, when there are (asymmetric) short loops,
we establish the existence of 
geodesics in $\Ham(M)$ by generalizing the nonsqueezing theorem to nontrivial 
 fibrations.

As in \S\ref{ss:geomint} consider a fibration $(P, \Om) \to B$ where $B = 
D$ or $S^2$.
We will say that {\bf the nonsqueezing theorem holds for} $(P,\Om)$ if
area$(P,\Om)$ constrains 
the radius of any embedded  symplectic ball $B^{2n+2}(r)$ 
 in $(P,\Om)$ by the inequality
$$
\pi r^{2} \;\le\; {\rm area\,} (P, \Om).
$$
For example, if $(P, \Om)$ is the product $(M\times D, \om 
+\al)$, where $\al$ is an area form on $D$, then $(P,\Om)$ has 
area equal to $\int_D\al$, so that this reduces to the usual 
nonsqueezing theorem.
 Similarly, we say that the nonsqueezing theorem holds
for the loop $\la$ if it holds for the corresponding fibration 
$(P_\la,\Om)\to S^2$
where $\Om$ is {\it any} $\om$-compatible symplectic form on $P_\la$.

The following result is proved in~\S\ref{sec:geod}. 

\begin{lemma}\label{le:nsgeo} \begin{itemize}
\item[(i)] Suppose that there is $\eps > 0$ such that
the nonsqueezing theorem holds for all loops 
$\la\in \pi_1(\Ham(M, \om))$ with $\trho_f([\la])\le \eps$, and let 
$(\tnu, \nu)$ be one of the pairs
$(\trho,\rho), (\trho_f,\rho_f)$.  Then there is a $C^2$-neighborhood 
$\Nn$ of $id$ in $\Ham(M, \om)$ such that every 
$\tnu$-minimizing path in $\Nn$ also minimizes $\nu$.

\item[(ii)]  The same statement holds for the pair
$(\tnu, \nu) = (\trho_f, \rho^+ + \rho^-)$ provided that 
the nonsqueezing theorem holds for all loops 
$\la$ with either $\ell^-(\la)$ or $\ell^+(\la) < \eps.$
\end{itemize} 
\end{lemma}

Using ideas from~\cite{LMP,Mcq} and 
Seidel~\cite{Seid2} we show in~\S\ref{ss:ns} that the hypothesis in (i)
above holds for all spherically rational symplectic manifolds $(M, \om)$.  
For such a manifold, the {\bf index of rationality} $q(M)$ is the smallest 
positive number $q$ such that $[\om](A)  \in q\Z$ for all $A\in 
H_2^S(M, \Z)$.
In the weakly exact case, we set $q(M) = \infty$.

\begin{prop}\label{prop:ns} If $(M, \om)$ is a spherically rational symplectic 
manifold with index of  rationality $q(M)$, the 
nonsqueezing theorem holds for all loops $\la$ in $\Ham(M, \om)$
with $\trho_f(\la) < q(M)/2.$
\end{prop}

For general manifolds $(M, \om)$ we establish the existence of 
$\rho$-geodesics by using a suitable modification of Lemma~\ref{le:nsgeo}. 
 Here the relevant quantity is the minimum size
of a class $B$ for which there is a nontrivial Gromov--Witten 
invariant $n_M(a,b,c; B)$.  Here nontrivial means that $B\ne 0$, and 
$a,b,c$ can be any elements in $H_{*}(M)$.  Thus we set
$$
\hbar = \hbar(M) = \min\{\om(B)> 0: \mbox{ some }\;n_M(a,b,c;B)\ne 0\},
$$
so that $\hbar = \infty$ if all nontrivial Gromov--Witten invariants 
vanish.  Note that $\hbar> 0$ for all $(M, \om)$: standard compactness results 
imply that for each $\ka$ there are only finitely many classes $B$ with
$\om(B) \le \ka$ that can be represented by a $J$-holomorphic curve 
for generic $J$, and it is only such classes that give rise to 
nonzero invariants.

\begin{prop}\label{prop:ns3}  Suppose that $\la$ is a loop in 
$\pi_1(\Ham(M), \om))$ such that $\ell^{\pm}(\la) < \hbar(M)/2$. Then there is 
 $\ka\in \R$ with $|\ka| \le \max(\ell^-(\la), \ell^+(\la))$
such that the radii of all symplectically 
embedded balls in $(P_{\pm\la}, \Om)$ are 
constrained by the inequalities
$$
\pi r^2 \le \area(P_{\la}, \Om) + \ka,\quad 
\pi r^2 \le \area(P_{-\la}, \Om) - \ka.
$$
In particular, if $\ell^{\pm}(\la) = 0$ then the nonsqueezing theorem holds 
for $\pm \la$.
\end{prop}

We call the above property {\bf weighted nonsqueezing}. 
For the proof see~\S\ref{ss:nsgeo}. 

\begin{cor}\label{cor:nsgeo}
Let $(\tnu, \nu)$ be one of the pairs
$(\trho,\rho), (\trho_f,\rho_f)$.
Suppose that the path $\ga$ has $\tnu$-length $< \hbar(M)/4$,
 minimizes $\tnu$ in $\THam(M)$ and is sufficiently short in the $C^2$-topology. 
 Then it is 
$\nu$-minimizing in $\Ham(M)$.
\end{cor}

The hypothesis that just one of  $\ell^+(\la), \ell^-(\la)$ is small
does not seem to give useful information towards proving the nonsqueezing 
theorem.  Therefore, to understand the $(\rho^+ + \rho^-)$-geodesics we use 
the following result.  Here we have written $c_1$ for the 
first Chern class of $(TM,J)$, where $J$ is any $\om$-tame almost complex 
structure on $M$.

\begin{prop}\label{prop:nsa} \begin{itemize}
\item[(i)] If $M$ is  weakly exact 
the nonsqueezing theorem holds for all loops $\la$ in $\Ham(M)$. 
\item[(ii)]  The same statement holds
if $c_1=0$ on $\pi_2(M)$ and all 
$3$-point Gromov--Witten 
invariants $n_{M}(a,b,c; B)$  vanish. 
\item[(iii)]
If $M = \CP^n$
the nonsqueezing theorem holds for 
the $(n+1)$st multiple $(n+1)[\la]$ of each loop $\la$. 
\end{itemize}
\end{prop}

These are sample results; our methods could doubtless be used to
find other manifolds for which 
the nonsqueezing theorem holds for all loops.  However, as is clear from 
Lemma~\ref{le:nsgeo} above, this is more than is needed to show 
that $(\rho^+ + \rho^-)$-geodesics exist.
There certainly are fibrations (such as the nontrivial $S^2$-bundle over $S^2$
with suitable $\Om$) 
for which the nonsqueezing theorem does not hold.  
Also  the nonsqueezing theorem may well fail for 
general fibrations $(P, \Om)\to D$ with nontrivial boundary monodromy.
(See Remark~\ref{rmk:nsthm} for further discussion on this point.)

\subsection{Calculating the seminorm $\rho\,^{+} + \rho\,^{-}$}\label{ss:exam}

This section describes the example mentioned in 
Theorem~\ref{thm:norm} (iii)
with $\rho(\phi)\ne\rho^{+}(\phi) + 
\rho^{-}(\phi)$.  In the case considered here,
 the path that gives the minimum of $\rho^{-}$ 
is not homotopic to the one that minimizes $\rho^{+}$.  Therefore, the 
example does not show that the norms $\trho$ and 
$\trho\,^{+}+\trho\,^{-}$
on $\THam(M)$ are different. (The latter question appears much more 
delicate: see \S\ref{ss:area}.)

\begin{prop}\label{prop:exam}  Suppose that 
 $\{\phi_{t}\}_{t\in [0,1]}$ is a loop $\la$ in $\Ham(M)$
 generated by a
function $H_{t}: M\to \R$ with $\Ll^{+}(H_{t}) \ne 
\Ll^{-}(H_{t}) $.  Suppose 
further that
$$
\Ll(\la) \;\;= \;\; r_{\trho}(M) = \inf \{\trho([\la]): \trho([\la]) > 
0, \;[\la]\in \pi_1(\Ham)\}.
$$
Then, if $\tau = \phi_{T}^{H}$ is the halfway point of this loop,
that is if $\Ll(\{H_{t}\}_{ t\in [0,T]}) = \Ll(\{H_{t}\}_{t\in [T,1]})$, 
$$
\rho(\tau) \;\;>\;\; \rho^{+}(\tau) + \rho^{-}(\tau).
$$
\end{prop}
\proof{}
 There are two natural paths to 
 $\tau$, namely $\be^{-}$ given by $\{\phi_{t}\}_{t\in [0,T]}$ and $\be^{+}$
 given by $\{\phi_{1-t}\}_{t\in [0,1-T]}$. 
 Without loss of 
 generality we may suppose  that $q =  \Ll^{-}(H_{t}) < \Ll^{+}(H_{t}) = p$.  
Hence
$$
q =  \Ll^{-}(H_{t}) \;=\; \Ll^{-}(\be^{-}) + \Ll^{+}(\be^{+})
\;\;<\;\;\Ll^{+}(\be^{-}) + \Ll^{-}(\be^{+}) = p.
$$
(This holds because the direction of $\be^{+}$ is the opposite of 
$\phi_{t}$.)
Thus
$$
\rho^{-}(\tau) +   \rho^{+}(\tau)\;\;\le \;\;\Ll^{-}(\be^{-}) + 
\Ll^{+}(\be^{+}) = q.
$$
On the other hand,  
 $\la$, by hypothesis, is an $\Ll$-minimizing representative of 
 its homotopy class.  This implies that both  paths $\be^{+}$ and $\be^{-}$ 
 minimize $\Ll$ in their homotopy class.  Further, since these 
 paths have length precisely $\Ll(\la)/2 = (p+q)/2$, any
 shorter path $\be'$  from $id$ to $\tau$  would create a nonconstant loop 
 $\la' = (-\be^+) * \be'$ in 
 $\Ham(M)$ with length $< \Ll(\la)$.  If $\ell_{\trho}([\la']) = 0$, 
 then we could alter $\be'$ by an arbitrarily short path to be 
 homotopic to $\be^{+}$ which contradicts the minimality of 
 $\Ll(\be^{+})$.  On the other hand, if $\ell_{\trho}([\la']) > 0$ 
 then we must have $\ell_{\trho}([\la']) \ge r_{\trho}(M) = \Ll(\la).$
 Hence the paths $\be^{\pm}$ must achieve the minimum of $\rho$, and 
 $$
 \rho(\tau)\; =\; (p+q)/2\; > \; q \;\ge\;  \rho^{-}(\tau) +  
 \rho^{+}(\tau).
 $$
\QED

In order to apply this argument to the norm $\rho_f$ instead of $\rho$  we need to 
start with a loop $\la$ whose length minimizes $\trho_f$ over $\pi_1(\Ham)$
and that reaches the halfway mark with respect to both $\Ll^-$ and $\Ll^+$ 
at the same time.  The latter condition is most easily achieved if $H_t$ is 
time independent, i.e. if $\la$ is given by a circle action.
Therefore, we have the following corollary.

\begin{cor}\label{cor:exam}  Suppose that a Hamiltonian $H$ with $\Ll^-(H) \ne \Ll^+(H)$
generates a circle action $\la = \{\phi_t^H\}$ and let $\tau = 
\phi_{1/2}^H$.  Then if
$\Ll(\la)$ minimizes $\trho_f$ over $\pi_1(\Ham)$ 
$$
\rho(\tau) \;\; = \;\;\rho_f(\tau)\;\;>\;\; \rho^{+}(\tau) + \rho^{-}(\tau).
$$
\end{cor}

The above results are easy.  Note that they imply that  the norms
$\rho$ (or $\rho_f$) and $\rho^{+} + \rho^{-}$ take 
different values on all elements
sufficiently close to $\tau$ in the $C^2$-topology. 

The next result is  harder.

\begin{prop}\label{prop:exam2}  There is a symplectic manifold
    $(M, \om)$ and a loop  $\la$ in $\Ham(M, \om)$
    that satisfies the conditions of Corollary~\ref{cor:exam}.
    \end{prop}
    
The main problem is to find a loop
 $\la$  whose length minimizes $\trho_f$.
It is shown in McDuff--Slimowitz~\cite{MSlim} 
that any semifree action of a circle achieves the minimum of $\trho$ in 
its homotopy class.  However, it is much more complicated to
estimate the lengths of all other loops. 
For this we will apply a 
method due to Seidel~\cite{Seid,Seid2} that uses 
the representation of $\pi_{1}(\Ham(M))$ 
on the quantum homology of $M$. 

There are few manifolds for which  
 $\pi_{1}(\Ham(M))$ is known. 
The easiest manifold to try would be $S^{2}$.  But this does not
work because the Hamiltonian 
for the generating loop of $\pi_1(\Ham(S^2)) = \pi_1(\SO(3))$ is symmetric.  
Therefore, following
 Polterovich~\cite{Pac},
 we will work out an explicit  example  on the
 one point blowup  $M_{*}$ of 
$\CP^{2}$.  We think of this as the region 
$$
\{(z_{1}, z_{2})\in \C^{2}: a^{2}\le |z_{1}|^{2} + |z_{2}|^{2}\le 1\}
$$
with boundaries collapsed along the Hopf flow.  
Abreu--McDuff~\cite{AM} show that $\pi_{1}(\Ham(M_{*}))=\Z$ for all 
$0 < a< 1$, with
 generator given by the rotation
$$
\al:\quad (z_{1}, z_{2})\mapsto (e^{-2\pi it}z_{1}, z_{2}), \quad 0\le t \le 1.
$$
However the $S^{1}$-action that satisfies the  hypothesis 
in Corollary~\ref{cor:exam} is
$$
\la:\quad (z_{1}, z_{2})\mapsto (e^{-2\pi it}z_{1}, e^{-2\pi it}z_{2}), \quad 0\le t \le 1.
$$
which is homotopic to $2\al$.  Thus 
$\tau(z_{1}, z_{2}) = (-z_{1}, -z_{2})$.  We will show 
that the hypotheses are satisfied for all $a$.
Details are  in~\S\ref{sec:exam}.

\MS\MS

\NI
{\bf Organization of the paper.\,\,}  \S2  proves  Proposition~\ref{prop:polt}
 about the 
geometric interpretations the the seminorms.  \S3 discusses geodesics.
The results here are based on the results in \S4 about the nonsqueezing 
theorem.  Finally, \S5 contains the calculations on the one point blow up 
of $\CP^2$ needed to prove Theorem~\ref{thm:norm}(iii).

\MS

\NI
{\bf Acknowledgements.\,\,}  I wish to thank Seidel for giving me 
a preliminary version of~\cite{Seid2}, Polterovich for some useful comments, 
Haydee Herrera and GuangCun Lu for pointing out some misprints in  an 
earlier version of this paper, and
Harvard University and the Courant Institute for 
providing a congenial atmosphere in which to  work.


\section{The geometric interpretation of $\rho\,^{\pm}$}\label{sec:geom}

Our arguments are based on the geometric approach to estimating
the Hofer norm  proposed in~\cite{LM}~II.  We will begin by reminding 
the reader of some definitions and notation from that paper.
 Throughout 
we will assume that $\om$ is normalized so 
that $\int_{M}\om^{n} = n!$

\subsection{The regions over and under the graph}\label{ss:gra}

In this section, unless explicit mention is made to the contrary
we will assume that $H_{t}$ is 
{\bf positively normalized}
in the following sense.  

We want to arrange that
the function 
$$
t\mapsto \min(t) = \min_{x \in M} H_t(x)
$$
is nonnegative and has arbitrarily small integral.
If the function $\min(t)$ is smooth, we simply  replace $H_{t}$ by
$H_{t} - \min(t)$.  This will be the case if $H_{t}$ has a fixed 
minimum.
In general, we replace $H_{t}$ by $H_{t} - m(t)$ where
$m(t)$ is a smooth function that is everywhere $\le \min(t)$
and is such that $\min(t) - m(t)$ has arbitrarily small integral.

Finally we reparametrize so that $H_{t}\equiv 0$
for $t$ near $0,1$.  It is easy to 
see that this can be done without changing $\Ll^{\pm}(H_{t})$.
In particular, every time independent Hamiltonian $H$ may be replaced by
one of the form $\be(t) H$ that satisfies the above condition 
and has the same length and time $1$-map as before.

We denote the graph $\Ga_H$ of $H_{t}$ by
$$
\Gamma_H = \{ ( x, t, H_t(x)) \} \subset M\times [0,1]  \times {\R} .
$$ 
For some small $\eps> 0$ choose a smooth function
$\mu'(t): [0,1] \rightarrow [-2\eps, 0]$ such that
$$
\int_0^1 \left(\min(t) - m(t) - \mu'(t)\right) dt = \eps.
$$
(This is possible provided that $m(t)$ is properly chosen.)
A {\bf thickening of the region under} $\Ga_H$ is 
$$
R_H^-(\eps) = \{ (x,t,h) \; \vline \; \mu'(t) \leq h \leq H_t(x) \} 
\subset M \times [0,1]  \times {\R}.
$$ 
Note that if $\mu'$ is suitably chosen so that its graph is
tangent to the lines $t= 0,1$ we may
arrange that $R_{H}^{-}(\eps)$ is a manifold with corners.
(Recall that $H_{t}\equiv 0$
for $t$ near $0,1$.)

 Similarly, we can define 
$R_H^+(\eps)$ to be a slight thickening of the region above 
$\Ga_{H}$: 
$$
R_H^+({\eps}) = \{ (x,t,h) \; \vline \; H_t(x) \leq h \leq \mu_H(t) \} 
\subset M  \times [0,1] \times \R
$$ 
where $\mu_H(t)$ is chosen so that   
$$
\mu_H(t) \geq \max_{x \in M} H_t(x) = \max(t),\qquad \int_0^1 (\mu_H(t) -
\max(t)) dt = \eps.
$$
  We define 
$$
R_H(2\eps)\; = \; R_H^-(\eps) \cup R_H^+(\eps)
\;  \subset 
M  \times [0,1] \times {\R},
$$ 
and equip $R_H^-(\eps)$, $R_H^+(\eps)$, and  
$R_H(2\eps)$ with the product symplectic 
form $\Omega_{0}= \omega + \al$ where $\al = dt \wedge dh$. 
In particular, for any Hamiltonian $H_{t}$, $(R_H(2\eps), \Omega_{0})$ 
is symplectomorphic to  the product
$(M \times D(a), \Om_{0})$ where 
$D(a)$ denotes the $2$-disc $D^2$ with area $a = \Ll(H) + 2\eps$. 

Clearly, there is a projection $\pi: 
R_H^\pm(\eps) \to D$  with fibers
of the form 
$$
\pi^{-1}(b) = \{(x, t_{b}, h_{b}(x)): x\in M, b\in D\}.
$$
Thus both spaces $(R_H^\pm(\eps), \Om_{0})$ fiber over the $2$-disc.
Further, if $X_t$ is the symplectic gradient of $H_t$, 
$$
\io(X_t + \p_t) (\om + dt\wedge dh) = -dH_t + dh
$$
vanishes on $\Ga_H$.  Therefore the monodromy along $\Ga_H$ in the 
direction of increasing $t$ is simply given by the flow $\phi_t^H$.
We can smooth the corners of the regions $R_{H}^{\pm}(\eps)$ 
so that the boundary monodromy is trivial except along $\Ga_H$.
Therefore, the boundary monodromy when taken in the direction 
corresponding to increasing $t$ is the element $\tphi$ in $\THam$.  
We have chosen here to orient the ambient space $M\times [0,1]\times \R$
in such a way that $t$ increases (resp. decreases)  
as one goes positively round the 
boundary of $R_H^+(\eps)$ (resp. $R_H^-(\eps)$).  (Here the boundary 
is oriented so that if one puts the outwards normal vector in front 
of a positively oriented basis for the boundary one gets a positively 
oriented basis for the ambient space.)
Hence the boundary monodromy of 
$(R_{H}^{+}(\eps), \Om_0)$ is $\tphi$ while that of
$(R_{H}^-(\eps), \Om_0)$ is $\tphi^{-1}$.

Finally, consider two Hamiltonians 
 $H_t$ and $K_t$ on $M$ such that  
$\phi_1^H = \phi_1^K$.   
There is a map $g:\Gamma_H \to \Gamma_K$ defined by  
$$
g(x,t,h) \;\;=\;\; \left(\phi_t^K \circ 
(\phi_t^H)^{-1}(x),\; t,\; h - H(x)+ K(\phi_t^K \circ 
(\phi_t^H)^{-1}(x))\right).
$$  
The above formula defines a symplectomorphism of $(M\times [0,1]\times 
\R, \Om_0)$ and so we use it to attach $R_H^-$ to $R_K^+$ to form the 
space\footnote
{
This definition is slightly different from that in~\cite{LM}~II to take 
account of the changed signs.} 
$$
(R_{K,H}(2\eps), \Omega_0) \;=\;
(R_K^+ (\eps)\cup_g R_H^-(\eps), \Om_0).
$$  
We assume that the functions $\mu'$ and $\mu_H$ are chosen so that 
$R_{K,H}(2\eps)$ is a smooth manifold with boundary. 
Note that $\Om_0$ has trivial monodromy round this boundary. 
Identifying this boundary to a point, we get a fibered space
$(P_{K,H}(2\eps), \Om_{0})\to S^{2}$.

In~\cite{LM}~II we only considered the case when
the loop
$\la= \phi_t^K \circ (\phi_t^H)^{-1}$ is contractible in $\Ham(M,\om)$,
in which case we showed that the fibration
$(P_{K,H}(2\eps), \Om_0)\to S^2$ is symplectically trivial.
(Such fibered spaces were called quasicylinders.)
In general, we identify $P_{K,H}(2\eps)$ with the space $P_{\la}$,
where $\la = \{\phi_t^K \circ (\phi_t^H)^{-1}(x)\}$ and 
so think of $\Om_{0}$ as a form on the fibration $P_{\la}\to S^2$.

\subsection{The area as norm}\label{ss:area}

We begin with a few words about fibrations of the form $(P,\Om)\to B$.
Even though we are only interested in the cases $B = D^2$ or $S^2$,
it is worthwhile to understand the geometric structure imposed by
the symplectic form $\Om$.  As pointed out earlier, $P\to B$ inherits the 
structure of a Hamiltonian fibration, i.e. its structural group may be 
reduced to $\Ham(M)$.  Conversely, it is shown in~\cite{LMh} and in
Ch. 6 of~\cite{MS}
that if a smooth fibration $\pi:P\to B$ with closed fiber 
$M$ has structural group 
$\Ham(M, \om)$ then the family of symplectic forms $\om_b$ on the fibers always 
has a closed extension.  
Moreover, if $B$ itself is a symplectic manifold, we may assume that this 
extension is itself a symplectic form $\Om$.
(The case $B = S^2$ is discussed 
in~\cite{Seid,LMP}.)
We claim:

\begin{quote}{\it the choice of such a form $\Om$ on 
a Hamiltonian fibration $\pi:P\to B$ with $B =  D$ or $S^2$ is 
equivalent to the choice of a pair $(\Ga, \al)$ where $\Ga$ is a
connection  on $\pi$ with Hamiltonian holonomy and $\al$ is an area form 
 on  $B$.}
\end{quote}  

To see this, first suppose that $M\to P\to B$ is 
 a  fibration with structural group 
$\Ham(M, \om)$.  One can always provide it with a connection
 $\Ga$  with Hamiltonian holonomy, 
and then use the  construction of 
Guillemin--Lerman--Sternberg to find a closed $2$-form 
$\tau_\Ga$ that extends the given forms on the fibers and induces the 
given connection $\Ga$.  
This form $\tau_\Ga$ is called the coupling form and is unique
(at least for closed bases) 
if one requires that the integral over the fiber of the class 
$[\tau_\Ga]^{n+1}$ is $0\in H^2(B)$, where $2n = \dim M$.
When $B= S^2$ this is equivalent to 
the condition 
\begin{equation}\label{eq:tau}
[\tau_\Ga]^{n+1} = 0\;\; \in  H^{2n+2}(P).
\end{equation}
(If $B = D$ uniqueness is given by  a relative version of this 
condition. For more detail see the proof of Lemma~\ref{le:polt3} below.)  
Further, because $M$ is compact one can obtain a symplectic 
extension $\Om_\Ga$ of the forms $\om_b$ on the fibers by adding to $\tau_\Ga$
 the pullback of a suitable area form on $B$.  
Hence the data $(\Ga, \al)$ does give rise to a symplectic extension 
$\Om = \Om_\Ga$ 
of the forms $\om_b$.   Note that $\Om_\Ga$ still 
 induces the connection $\Ga$  on $P\to B$.
Hence, if we start from $\Om$ and let $\Ga$ be the corresponding 
connection, both  $\Om$ and the coupling form $\tau_\Ga$ 
induce the same connection.   Since the uniqueness condition for 
$\tau_\Ga$  is cohomological, it is satisfied by some form of the type
$\Om - \pi^*(\al)$.  Hence $\Om$  may be written as
$$
\Om = \tau_\Ga + \pi^*(\al)
$$
where $\tau\,_\Ga$ is the coupling form for the connection $\Ga$ defined 
by $\Om$ and $\al$ is an area form on $B$.  This proves the claim.
(For more details, see for example~\cite{Pk,MS,LMh}.)

\begin{lemma} If $\Om = \tau_\Ga + \pi^*(\al)$ as above, then
$$
{\rm area}\,(P, \Om) = \int_B \al.
$$
\end{lemma}
\proof{}  This follows immediately from
 condition~(\ref{eq:tau}) above and the normalization condition
 $\int \om^{n} = n!$\QED

We next explain
Polterovich's ideas from~\cite{Peg,Pbk} that give geometric 
interpretations of the seminorms $\trho^{\pm}$.  Recall from 
Definition~\ref{def:area} that
$$
\ta^{+}(\tphi) = \inf{\area (P,\Om)}, \qquad \ta^{-}(\tphi) = 
a^{+}(\tphi^{{-1}}),
$$
where the infimum is taken over all $\om$-compatible fibered spaces
$(P, \Om) \to D$ with boundary monodromy equal to $\tphi$.  
A simple calculation shows that
$$
\area\,(R_H^+(\eps), \Om) = \Ll^{+}(H_{t}) + \eps,\quad
\area\,(R_H^-(\eps), \Om) = \Ll^{-}(H_{t}) + \eps.
$$
Therefore,
$$
\ta^{+}(\tphi) \;\le \;\trho\,^{+}(\tphi),\qquad
\ta^{-}(\tphi) \;\le\; \trho\,^{-}(\tphi). 
$$
The following result combines Polterovich~\cite{Pk}\S3.3 and~\cite{Peg}\S3.3.
We give a proof here partly because he considers loops rather than 
paths and partly because we  simplify his argument by
avoiding the use of  K-area.

\begin{lemma}\label{le:polt3}
$\ta^{+}(\tphi) = \trho\,^{+}(\tphi)$ and
$\ta^{-}(\tphi) = \trho\,^{-}(\tphi)$. 
\end{lemma} 
\proof{}  By symmetry it suffices to prove the former 
result, which will follow if we  show that 
$$
\ta^{+}(\tphi) \;\;\ge\;\; \trho\,^{+}(\tphi).
$$
Suppose to the contrary that we are given a fibration $(P,\Om)\to D$ 
with area $ < \trho\,^{+}(\tphi)$ and monodromy $\tphi$.
By Moser's theorem we may isotop $\Om$ so that it is a product 
in some neighborhood $\pi^{-1}(\Nn)$ of the
base fiber $M_*$. Identify the 
base $D$ with the unit  square $K = \{0\le x,y \le 1\}$ taking 
$\Nn$ to a neighborhood of $\p' K = \p K - \{1\}\times (0,1)$, 
and then identify $P$ with 
$K\times M$ by parallel translating along the lines $\{(x,y): x\in [0,1]\}$. 
In these coordinates, the form $\Om$ may be written as 
$$
\Om\;\; = \;\; \om +  d_MF'\wedge dy + L' dx\wedge dy
$$
where $F', L'$ are suitable functions on $K\times M$ and $d_M$ 
denotes the fiberwise exterior derivative.  
Because $\Om$ is a product near  $\pi^{-1}(\p' K)$, 
$d_MF' = 0$ there and $L'$ reduces to a function of $x,y$ only. By subtracting 
a suitable function $c(x,y)$ from $F'$  we can arrange that $F = F' - c(x,y)$ 
has zero mean
 on each fiber $\pi^{-1}(x,y)$ and then write $L' + \p_xc(x,y)$ as $-L + 
 a(x,y)$ where 
$L$ also has zero fiberwise means.  Thus 
\begin{equation}\label{eq:form}
\Om\;\; = \;\; \om +  d_MF\wedge dy - L dx\wedge dy + a(x,y)dx\wedge dy,
\end{equation}
where both $F$ and $L$ vanish near $\pi^{-1}(\p' K)$ and have zero 
fiberwise means.
 Since $\Om$ is symplectic it must be positive on the
 $2$-dimensional distribution \,\,{\it Hor}.  Hence we must have $
-L(x,y,z) + a(x,y)> 0$ for all $x,y\in K, z\in M$.
Moreover, it is easy to see that 
$\area (P, \Om) = \int a(x,y) dx \wedge dy.$  Hence
\begin{equation}\label{eq:a2}
\int \max_{z\in M} L(x,y,z) \,dx\wedge dy < \int a(x,y)\,dx\wedge dy = 
\area(P,\Om).
\end{equation}

We claim that $-L$ is the curvature of the induced
connection $\Om_\Ga$.  To see this, consider the vector fields
$X = \p_x, Y = \p_y - \sgrad\, F$ on $P$ that are
 the  horizontal lifts of $\p_x, \p_y$.\footnote
{
Here the symplectic gradient  $\sgrad\, F$ is defined by setting
  $\om(\sgrad \,F, \cdot) = -d_MF(\cdot)$.}
It is easy to check that their commutator $[X,Y] = XY - YX$ is 
vertical and that
$$
[X,Y] = -\sgrad(\p_xF) = \sgrad\, L 
$$
on each fiber $\pi^{-1}(x,y)$ as claimed.  (In fact, the first three terms 
in~(\ref{eq:form})  make up the coupling form $\tau_{\Ga}$.)

Now let $f_{s}\in \Ham(M)$ be the monodomy of $\Om_\Ga$ along the path
$t\mapsto (s,t), t\in [0,1]$.  (This is well defined because all  
fibers have a natural identification with $M$.)
The path $s\mapsto f_{s}$ is a Hamiltonian isotopy
from the identity to $\phi = f_{1}$, and it is easy to see that it is 
homotopic to the original path $\tphi$ 
given by parallel transport along $t\mapsto 
[1,t]$.  (As an intermediate path in the homotopy take
the monodromies along $t\mapsto (s,t), t\in [0,T]$ for $s\in [0,1]$
followed by the lift
of $\tphi_t, t\in [T,1]$, to $\THam$.)  Therefore 
 $\Ll^+(f_s) \ge \trho^+(\tphi)$, and we will derive a contradiction 
 by estimating $\Ll^+(f_s) $.

To this end, let $X^s, Y^t$ be the (partially defined) flows of the vector 
fields $X,Y$ on $P$ and set $h_{s,t} = Y^t X^s$.  Consider the 
$2$-parameter family of (partially defined)
vector fields $v_{s,t}$ on $P$ given by
$$
v_{s,t} = \p_s h_{s,t} = Y^t_*(X)\quad\mbox{on }\;\; {\rm Im\,}h_{s,t}.
$$
In particular $v_{s,1}(x,y)$ is defined when $y=1, s\le x$.
Since $f_s = h_{s,1}$ we are interested in calculating the 
vertical part of $v_{s,1}(s,1,z)$.  Since the points 
with $y = 1$ are in ${\rm Im\,}h_{s,t}$ for all $(s,t)$ we may write
\begin{eqnarray*}
v_{s,1} & = & \int_0^1 \p_t(v_{s,t})\; dt + v_{s,0}\\
& = & \int_0^1 Y^t_*([X,Y])\; dt + \p_x.
\end{eqnarray*}
We saw above that $[X,Y] = \sgrad\, L$. Hence $Y^t_*([Y,X]) = \sgrad(L\circ 
(Y^t)^{-1})$ and
\begin{eqnarray*}
v_{s,1}(s,1,z) & = & \int_0^1 \sgrad(L((Y^t)^{-1}(s,1,z)) \; dt + \p_x\\
& = & \sgrad \int_0^1\,L(s,1-t, (Y_v^t)^{-1}(z))\;dt + \p_x
\end{eqnarray*}
where $Y_v^t$ denotes the vertical part of $Y^t$.
Hence the  Hamiltonian $H_{s}$
that  generates the path $f_{s}, s\in [0,1],$ and has zero mean
satisfies the inequality
$$
H_{s} (z)\le 
\int \left(\max_{z\in  M} L(s,t,z)\right) dt < \int a(s,t) dt,
$$
since $L(s,t,z)< a(s,t)$ by~(\ref{eq:a2}).
Thus  $\trho\,^+(\tphi)\le \area P$, contrary to hypothesis.\QED

Using this, we can interpret the seminorms $\rho^{+} + \rho^{-}$ and 
$\rho_{f}$ in terms of the fibrations
$P_{K,H}(\eps)\to S^{2}$ considered at the end of~\S\ref{ss:gra}.

\begin{cor} \label{cor:polt} \begin{itemize}
\item[(i)]  $\rho^{+}(\phi) + \rho^{-}(\phi)$ is the minimum 
of the areas of the fibrations
$(P_{K,H}(\eps), \Om_{0})\to S^{2}$ taken over all $\eps>0$ and 
all pairs  $(H_{t}, K_{t})$
of Hamiltonians
with time $1$ map $\phi$.
\item[(ii)]  $\rho_{f}(\phi)$ is the minimum of these areas over the 
set of pairs that define homotopic flows.  
\end{itemize}
\end{cor}

The proof of   Proposition~\ref{prop:polt} is now immediate.\MS

Now let us consider the geometric interpretation of the Hofer norm
and compare it with that for $\rho_{f}$.  The above corollary implies 
that $\rho_{f}$ is the minimum of the areas of the symplectically 
trivial fibrations 
$(R_{K,H}, \Om_{0})$ that contain a hypersurface $\Ga_{H} = \Ga_{K}$ with 
monodromy $\phi$.  The sets $R_{K,H}$ are constructed as subsets of $(M\times 
\R^{2},\om + dt\wedge dh)$ that have trivial monodromy round the boundary.  This 
means that the flow lines of the characteristic flow round the 
boundary are all closed. However, they are not constant in the $M$ 
direction, i.e. with respect to the given trivialization of the fibers they go 
round the loops  $t\mapsto \phi_{t}^{K}\circ (\phi_{t}^{H})^{{-1}}(x)$.
The Hofer norm, on the other hand, is the minimum area of
the cylinders $R_{H}(2\eps)$ which sit inside $(M\times \R^{2}, \om + 
dt\wedge dh)$
as product regions of the form $M\times U_{H}$.  Hence the monodromy
of its boundary is constant with respect to the given trivialization.
Moreover the 
hypersurface $\Ga_{H}$ in $R_{H}(2\eps)$ that has monodromy $\phi$ is 
a graph over its front and back faces (defined by $h = \mu_{H}(t)$,
$h = \mu'(t)$).  For general $M$ 
is it not at all clear that these are the same minimizing sets.
Even when $M = S^{2}$ and we know that $(R_{K,H}, \Om_{0})$ is 
symplectomorphic to a product subset of $(M\times \R^{2}, 
\om + dt\wedge dh)$,
it is not obvious that the resulting hypersurfaces in this subset can be 
assumed to be graphs over the front and back faces.

\MS

\NI
{\bf Proof of Lemma~\ref{le:polt2}}.\,\,
By Proposition~\ref{prop:polt} $\trho^+([\la])$ is measured by minimizing
the area of fibrations 
 $(P, \Om)\to D$  with boundary 
monodromy in the class $[\la]$.  Therefore, the result will follow if we 
set up a correspondence between this minimizing set and
the set of fibrations $(P_\la, \Om) \to S^2$ with clutching loop $[\la]$.

Given  $(P, \Om)\to D$  with monodromy $[\la]$ there is an associated
 fibration $Q\to S^2$ defined by 
identifying the boundary $\p P$ to a single fiber via the characteristic 
flow.  Moreover, since the restriction of $\Om$ to $\p P$ determines 
$\Om$ near $\p P$, $\Om$ descends to $Q$.  We therefore get a 
fibration $(Q, \Om) \to S^2 = D/\p D$, with fiber over the base point
$* = \p D$ identified to $M$.  It remains to observe that the 
clutching function of this 
fibration $Q\to S^2$ is   $[\la]\in \pi_1(\Ham)$. 
Conversely, given such a fibration $(Q, \Om)\to S^2$ identify $(S^2, *)$ with 
$(D/\p D, \p D)$ and consider the pullback of $Q$ to $D$.  This will be a 
fibration over $D$ with boundary monodromy in the class $[\la]$, as required.
 \QED

\section{Geodesics}\label{sec:geod}

This section contains proofs of the results in \S\ref{ss:geod} assuming 
the results stated in \S\ref{ss:sloop} about the nonsqueezing theorem.

Let $\{\phi_{t}\}_{t\in [0,1]}$ be a path 
 generated by a Hamiltonian $H_{t}$
that is positively normalized as in \S\ref{ss:gra}.  The following 
easy result was proved  in~\cite{LM}~II~Lemma~3.2.
Recall that the capacity of a ball of radius $r$ is $\pi r^{2}$.

\begin{lemma}\label{le:1}
If $H_{t}$ is sufficiently small in the $C^{2}$-norm and has a fixed 
maximum (resp. minimum), then, for all 
$\eps > 0$ it is 
possible to embed a ball of capacity $\Ll(H_{t})$ in 
$R_{H}^{-}(\eps)$ (resp. $R_{H}^{+}(\eps)$).
\end{lemma}

\begin{lemma}\label{le:2}  Suppose that a ball of capacity $\Ll(H_{t})$ can be 
embedded in  $R_{H}^{-}(\eps)$ for all $\eps$.  Then:
\begin{itemize}
\item[(i)]  the 
corresponding path $\tphi_{t}= \phi_{t}^{H}$  in $\THam$ minimizes
$\trho\,^{+}(\tphi)$;
\item[(ii)] the path $ \phi_{t}^{H}$ also minimizes $\rho^{+}(\phi)$
provided that the nonsqueezing theorem holds for all
 loops $\la\in \pi_{1}(\Ham).$
\end{itemize}
\end{lemma}
\proof{}  Suppose that $K_{t}$ generates a homotopic path $\psi_{t}$ with 
$\psi_{1} = \phi_{1}$  and that $\Ll^{+}(K_{t}) < \Ll^{+}(H_{t})$.
Then,  the corresponding fibered space $(P_{K,H}(2\eps), \Om) \to S^{2}$ has 
area $< \Ll(H_{t}) =  \Ll^-(H_{t}) +  \Ll^+(H_{t}) $ provided that $\eps$ is
sufficiently small.  On the other hand it contains a ball of capacity 
$\Ll(H_{t})$.  But this contradicts the 
nonsqueezing theorem for the symplectically
trivial bundle $(P_{K,H}(2\eps), \Om) \to S^{2}$.
This proves (i).  (ii) is also immediate, because the hypothesis means 
that we can apply the nonsqueezing theorem to any bundle of the form
 $(P_{K,H}(2\eps), \Om) \to S^{2}$.
\QED

We now show that paths that minimize $\trho\,^{+}$ must have fixed 
maxima.  

\begin{lemma}\label{le:3}  Suppose that $H_{t}$ 
 does not have a fixed maximum.
    Then  the 
corresponding path $\phi_{t}= \phi_{t}^{H}$  in $\THam$ does not 
minimize $\trho\,^{+}$.
\end{lemma}
\proof{}  The analogous result for the two sided norm $\rho$ 
was proved in~\cite{LM}~I~Proposition~2.1 by a simple curve shortening 
procedure. There we constructed a perturbation $\Psi_{t}^{\eps}$ such that:
\MS

\NI$\bullet$  $\Psi_{0}^{\eps} = \Psi_{1}^{\eps} = id.$
\smallskip

\NI$\bullet$  if $\phi_{t}^{\eps} = \Psi_{t}^{\eps}\circ \phi_{t}$ for 
all $t$, then the generating Hamiltonian $H_{t}^{\eps}$ for
$\phi_{t}^{\eps}$ satisfies the conditions:
$$
\min_{x\in M} (H_{t}^{\eps}) = \min_{x\in M} (H_{t}), \qquad 
\max_{x\in M} (H_{t}^{\eps}) \le \max_{x\in M} (H_{t}),
$$
where strict inequality holds for the maximum on some interval $|t - 
t_{0}| < \eps$. 
\MS

This implies that the two sided 
length $\Ll$  is smaller on $\phi_{t}^{\eps}$. To ensure that
$
\Ll^{+}(\phi_{t}^{\eps}) < \Ll^{+}(\phi_{t})$ it suffices to arrange 
in addition that
$$
\int H_{t}^{\eps}\, \om^{n} = \int H_{t}\, \om^{n} = 0,\qquad t \in 
[0,1].
$$
Clearly we can arrange that $H_{t}$ has zero mean, and 
$H_{t}^{\eps}$ will have too as long as we choose the functions 
$K_{j}$ that generate the perturbation $\Psi_{t}^{\eps}$
to have zero mean.  

The $K_j$ were constructed as follows.
We chose $t_0 < t_1 < \dots < t_k$ so that
$\cap_j X_j = \emptyset$, where $X_j = {\rm maxset\,}H_{t_j} \ne M$, and then 
chose a partition of unity $\{\be_j\}$ subordinate to the cover 
$M - \Nn_\ka(X_0), M-X_1, \dots, M- X_k$, where $\Nn_\ka(X_0)$ is the 
$\ka$-neighborhood of $X_0$.  Fix $\de > 0$ so that 
$$
\Nn_{\ka/2} (X_0)\;\;\subset\;\;
\cup_{j\ge 1}\be_j^{-1}([\de,1]),
$$
 and choose 
functions $K_j\le 0$ with support in $\be_j^{-1}([\de/2,1])$ that are constant 
and $< 0$ on $\be_j^{-1}([\de,1])$.   
Then define $\Psi_t^\eps$ as a smoothing of the 
following path, where  the flow of 
$K_j$ is denoted $\psi_t^j$:
\smallskip

\NI$\bullet$
$\Psi_t^\eps = id$ for $t < t_0 - \eps$ and then equals 
$\psi_s^1\circ\dots\circ\psi_s^k$ for $0 \le s= t-t_0 + \eps \le 2\eps$;
\smallskip

\NI$\bullet$
$\Psi_t^\eps$ is constant when $|t - t_j| > \eps$;
\smallskip

\NI$\bullet$
when $|t-t_j|\le \eps$ for some $j> 0$, $\Psi_t^\eps$ has the form 
$$
(\psi_s^j)^{-1}\Psi_{t_j - \eps},\quad\mbox{where }\; s = t - t_j+\eps.
$$ 

It is not hard to see that this satisfies all the requirements for small 
enough $\eps$. In particular, because $\sum_{j\ge1} K_j(x) < 0$ for 
$x\in \Nn_\ka(X_0)$, $\max H_t^\eps < \max H_t$ for $|t-t_0|< \eps$.
(Details are in~\cite{LM}~I.)

 In the present situation, we need to allow each $K_j$ to be 
positive somewhere  so that it can have zero mean.  
If  $M - \Nn_{\ka/2}(X_j \cup X_0)$ is nonempty for small $\ka$,
then there is no problem;
the argument goes through as before if $K_j$ is small and positive 
on such a set.   
Conceivably, we have to choose some of the $X_j$
 so that $M = X_0 \cup X_j$ in order to achieve that 
$\cap_{j=0}^k X_j = \emptyset$.
In this case the frontier of $X_0$ lies entirely in $X_j$.  So there must 
be an open set $U_j\subset X_0\cap X_j$ that lies outside some other $X_k$.  
So for each such $j$ we 
let  $K_j$ be positive in this open set $U_j$ and make $K_k$ sufficiently 
negative in $U_j$ to compensate.\QED

\begin{cor}  Proposition~\ref{prop:geod} holds.
\end{cor}
\proof{}
It follows from Lemmas~\ref{le:1}, \ref{le:2} above that a sufficiently 
short piece $\ga$ of any path with a fixed maximum (resp.  
minimum) 
minimizes $\trho\,^{+}$ (resp. $\trho\,^{-}$).  Conversely, if a path 
does not have a fixed maximum it cannot minimize  $\trho\,^{+}$.
\QED
    
\NI
{\bf Proof of Lemma~\ref{le:nsgeo}}.\,\,
To prove (i) we have to show that there is some $C^2$-neighborhood of $id$
in $\Ham$ such that every $\tnu$-minimizing path $\ga$  in $\Nn$ actually 
minimizes $\nu$, where $(\tnu, \nu)$ is 
$(\trho,\rho)$ or $(\trho_f,\rho_f)$.  
Clearly, it suffices to prove this for paths that start at $id$.

In both cases, we know from 
Lemma~\ref{le:3} above that  $\ga = \{\phi_t^H\}$ 
 has fixed extrema.  Therefore, by 
Lemma~\ref{le:1} we may choose $\Nn$ so that 
a ball of capacity $\Ll(\ga)$ embeds in $R_H^\pm(\de)$ 
for all $\de > 0$.
We may also choose $\Nn$ so that $\Ll(\ga) = \trho(\ga) < \eps/2$.
We claim that such $\ga$ must minimize both $\rho$ and $\rho_{f}$.  
We will carry out the argument for $\rho_{f}$ since it is slightly 
more complicated.

Suppose that $\ga$ does not minimize $\rho_{f}$.
Then there are paths 
 $\psi_{t}^{\pm}$ from $id$ to $\phi = \phi_1^H$  
generated by $K_{t}^{\pm}$ such that
 $$
\Ll^{+}(K_{t}^{+}) + \Ll^{-}(K_{t}^{-}) \;<\; \Ll(\ga) \;=\; 
\Ll^{+}(H_{t}) + \Ll^{-}(H_{t}).
$$ 
Therefore at least one of the following inequalities must hold:
$$
\Ll^{+}(K_{t}^{+}) \;\; < \;\; \Ll^{+}(H_{t}),
\qquad \Ll^{-}(K_{t}^{-}) \;\; < \;\; \Ll^{-}(H_{t}),
$$
say the former.   
Note that the two functions $K_{t}^{+}$ and $K_{t}^{-}$ may be 
different.  (In the case of $\rho$ they would be the same.)
However, because we are dealing with the norm $\rho_{f}$
rather than $\rho^{+} + \rho^{-}$ the paths $\be = \be^\pm$
that  they generate are homotopic.
Hence the  fibrations $P_{K^{+},H}, P_{H,K^{-}}$ correspond to loops $\la, 
-\la$ that are mutual inverses, where $\la = \be * (-\ga)$.
  Note also that for suitably small $\de$
$$
\trho^+([\la])\;\; \le \;\; \area(P_{K^{+},H}(\de), \Om_{0})\;\; < \;\;\Ll(\ga) 
\;\;<\;\; 
\eps/2,
$$
while
$$
\begin{array}{lcl}
\trho_f([\la]) & = & \trho^+([\la])  + \trho^-([\la]) \\
& \le &
\area(P_{K^{+},H}(\de), \Om_0) + \area(P_{H,K^{-}}(\de), \Om_{0}) \\
 & =  & \Ll(H) + \Ll^+(K_t^+) + \Ll^-(K_t^-) + 2\de\\
 &\le & \eps.
\end{array}
$$
Therefore by hypothesis the nonsqueezing theorem 
 holds for $\pm \la$.  Hence $(P_{K^{+},H}(\de), \Om_{0})$ cannot contain a ball 
 of capacity $\Ll(\ga)$.    This contradiction shows that $\ga$ must minimize 
 $\rho_{f}$.  
 Hence (i) holds. The proof of (ii) is similar: compare Lemma~\ref{le:asloop}.
\QED

\NI
{\bf Proof of Theorem~\ref{thm:geod}.}\,\,
When $(M, \om)$ is spherically rational
 (i) follows from Proposition~\ref{prop:geod}
by  Lemma~\ref{le:nsgeo} and Proposition~\ref{prop:ns}.  For the general case see~\S\ref{ss:nsgeo}.
The statement in (ii) about weakly exact $M$ also follows by a 
similar argument since, by Proposition~\ref{prop:nsa},
 the nonsqueezing theorem now holds 
for all loops.  
Therefore, it remains to consider the case $M = \CP^{n}$.

Let $\{\phi_{t}\}_{t\in [0,1]}$ be a path in $\Ham(\CP^{n})$ with 
a fixed maximum and minimum.  We must show that
 sufficiently short pieces of it 
 minimize $\rho^{+} + \rho^{-}$. In the following we denote by
 $H_t * K_t$  the 
 Hamiltonian $H_t + K_t\circ (\phi_t^H)^{-1}$ that generates the 
 composite $\phi_t^H\circ \phi_t^K$.  Further we write $m = n+1$.

First choose $\eps$ so that each piece 
$$
\ga_{a} = \{ 
\phi_{a+t} \phi_{a}^{-1}\}_{t \in [a, a+\eps]}
$$
is generated by a Hamiltonian $H_{t}^{a}$ for which the $m$-fold 
composite
$$
(H_t^a)^{*m} = H_{t}^{a}*\dots * H_t^a
$$
satisfies the conditions of Lemma~\ref{le:1}. 
Since $H_t^a$ has fixed extrema, 
$$
\Ll((H_t^a)^{*m}) = m\Ll(\ga_{a}).
$$
Therefore for each 
$\ga_{a}$ a ball of capacity $m\Ll(\ga_{a})$ can be embedded
in $R_{(H^a)^{*m}}^{\pm}$.
If some such piece $\ga_{a}$ does not minimize $\Ll^{+}(\ga_{a})$, 
then there is a shorter path $\ga'$ 
with the same endpoints generated by some $K_{t}$.
Let $\la$ be the loop $\ga'* (-\ga_{a})$ and consider the 
fibration
$$
(P_{m\la}, \Om)\;\; =\;\; (P_{(K^{*m}, (H^{a})^{*m})}, 
\Om_{0})\;\;\to\;\; S^{2}.
$$
The nonsqueezing theorem holds for this fibration by 
Proposition~\ref{prop:nsa}.  On the other hand it has area strictly 
less than   $m\Ll(\ga_{a})$ while containing embedded balls of 
this capacity, a contradiction. A similar argument applies if
$\Ll^-$ is not minimized.  Hence result.\QED
\MS

\NI
{\bf Proof of Proposition~\ref{prop:nbhd}}\,\,
This states that there is a $C^{2}$-neighborhood $\Nn$ of $id$
on which the two norms $\rho_{f}$ 
and $\rho$ agree.   It is shown in~\cite{LM}~II Proposition~5.11
that we may suppose that every $\phi\in \Nn$ is generated by some 
Hamiltonian $H_{t}$  with a fixed maximum and minimum.   Moreover, the 
length of this path is given by a Banach space norm.
If in addition $H_{t}$ is sufficiently $C^{2}$-small we may apply the 
arguments above to conclude that the path minimizes both $\rho$ and 
$\rho_{f}$.\QED

\begin{rmk}\rm
Here is an explicit description for 
the size of this 
neighborhood $\Nn$.  First,  to get an explicit 
choice of generating Hamiltonian $H_t$ for 
each $\phi\in \Nn$,  choose an identification of
a neighborhood of the diagonal in $M\times M$ with the zero section in 
$T^*M$.  Secondly,  to get uniform bounds
for the embedded balls 
of Lemma~\ref{le:1}, we can argue as follows.  Choose for some $\eps > 0$  
a family of symplectic embeddings
$\io_x:B^{2n}(\eps)\to M$ such that $\io_x(0) = x$ for all $x\in M$.  
This family should 
have uniformly bounded second derivatives, though it need not depend
 smoothly (or even continuously) on $x$.  It can be constructed from a 
 finite covering of $M$ by Darboux charts $B^{2n}(r_i)\to M$ chosen so that
the images of the subballs $B^{2n}(r_i - \eps)$ also cover $M$.
 Then, we require $\Ll(H_t)\le \eps$.  Further,
for each $H_t$ there must be at least one fixed maximum $x^+$ 
and one fixed minimum $x^-$ such that
 the functions $H_t\circ\io_{x^\pm}$ have sufficiently 
small  $C^2$-norm for the embedding techniques of~\cite{LM}~II to work.
 Finally, we need $\Ll(H_t)< \hbar/4$ in order for 
Corollary~\ref{cor:nsgeo} to hold. 
 \end{rmk}

\NI
{\bf Proof of Proposition~\ref{prop:good}.\,\,}
The first statement is proved in~\cite{MSlim}.  Although that paper 
only mentions the usual Hofer norm, the capacity--area inequality 
proved there shows that $\{\phi_{t}^{H}\}_{0\le t\le 1}$ 
minimizes both $\trho\,^{+}$
and $\trho\,^{-}$.  (Argue as in Lemma~\ref{le:Jsph} above using
~\cite{MSlim} Proposition 2.4.)

To prove the next statement it suffices
to show that if $M$ is weakly exact then all the norms 
$\rho^+ + \rho^-$, $\rho_f$ and $\rho$ are minimized as well.
To do this we simply have to show that the arguments in~\cite{MSlim} 
apply when the fibration $(P_{K,H}, \Om_{0})\to S^{2}$
is not symplectically trivial.  However all we used about this 
fibration is that there is a section $\si_{A}$ with 
$u_{\la}(\si_{A})\le 0$ such that the 
Gromov--Witten invariant $n_{P}([M], [M], pt;\si_{A})$ is nonzero:
see~\cite{MSlim}~\S3.1.  Since  this holds for weakly exact $M$ by 
Lemma~\ref{le:Jsph}, the 
result follows.\QED

\begin{rmk}\rm\begin{itemize}\item[(i)]
 The above argument shows that Proposition~\ref{prop:good}
holds for every $M$ for which 
 all loops $\la$ have good sections of positive weight 
in the sense of 
 Definition~\ref{def:sns} below.  It was based on~\cite{MSlim} but one could 
 equally use Entov~\cite{En}.
 \item[(ii)]  The proof of Theorem~\ref{thm:geod}(ii) for $\CP^{n}$ 
 also applies to any manifold $M$ that is not spherically 
 monotone\footnote
 {
 $(M,\om)$ is said to be spherically monotone
 if there is $\ka$ such that  $c_1-\ka [\om]$ vanishes on all classes 
 in $\pi_2(M).$ 
 }
 and such that
 $\pi_{1}(\Ham(M))$  is finite, or more generally has finite image 
 under the representation $\Psi$ defined below.
 See Remark~\ref{rmk:nsthm}(i).
 \end{itemize}
\end{rmk}


\section{Quantum homology and the nonsqueezing theorem}\label{sec:sloop} 
 
In~\cite{Mcq} Proposition~1.1 we established the nonsqueezing theorem for
some fibrations $P_{\ga}\to S^{2}$ by using the modified Seidel representation
$$
\Psi: \pi_{1}(\Ham(M)) \to (QH_{\ev}(M, \La))^{\times }.
$$
Here $(QH_{\ev}(M, \La))^{\times}$ denotes the 
group of units in the even part of the small quantum homology of $M$
with coefficients in a real Novikov ring $\La$. 
Since we use the same approach in this paper, we will begin by
recalling some definitions from~\cite{LMP,Mcq}.

\subsection{Preliminaries}

Set $c_1 = c_1(TM)\in H^2(M,\Z)$.  Let $\La$ be the   Novikov 
ring of
the group $\Hh = H_2^S(M,\R)/\!\!\sim$  with valuation $I_\om$ where 
$B\sim B'$
if $\om(B-B') = c_1(B-B') = 0$.
 Thus $\La$ is the completion of the rational group ring
 \footnote
 {
In~\cite{Mcq,LMP} we distinguished 
between the integral version of $\La$ 
which is generated by the integral elements of $\Hh$ 
and the real Novikov ring that we are now calling $\La$.  It is not 
necessary to do that here.}
 of $\Hh$ with elements of the form 
$$
\sum_{B\in \Hh} q_B\; e^B
$$
where for each $\ka$ there are only finitely many nonzero
$q_B\in \Q$ with $\om(B) > - \ka$.
  Set 
$$
QH_{*}(M) = QH_*(M,\La)
= H_*(M)\otimes\La.
$$
We may define an $\R$ grading on $QH_*(M,\La)$ by setting
$$
\deg(a\otimes e^B) = \deg(a) + 2c_1(B),
$$
and can also think of
$QH_*(M,\La)$ as $\Z/2\Z$-graded with 
$$
QH_{\ev} =  
H_{\ev}(M)\otimes\La, \quad QH_{odd} =  
H_{odd}(M)\otimes\La.
$$

Recall that the quantum intersection product 
$$
a*_Mb\in QH_{i+j - 2n}(M), \qquad \mbox {for }\; a\in H_i(M), b\in 
H_j(M)
$$
 is defined as follows:
\begin{equation}\label{eq:qm0}
a*_Mb = \sum_{B\in \Hh} (a*_Mb)_B\otimes e^{-B},
\end{equation} 
where  $(a*_Mb)_B\in H_{i+j- 2n+2c_1(B)}(M)$ is defined by the 
requirement
that 
\begin{equation}\label{eq:qm}
(a*_Mb)_B\,\cdot_M\, c = n_M(a,b,c;B) \quad\mbox{ for all }\;c\in 
H_*(M).
\end{equation}
Here  $n_{M}(a,b,c; B)$ denotes the Gromov--Witten invariant that counts 
the number of $B$-spheres in $M$ meeting the cycles $a,b,c\in 
H_{*}(M)$,
and we have written $\cdot_M$ for the usual  intersection
pairing on $H_*(M) = H_*(M, \Q)$.  Thus
$a\cdot_M b = 0$ unless $\dim(a) + \dim (b) = 2n$ in which case it is 
the
algebraic number of intersection points of the cycles. The product  
$*_M$ is
extended to $QH_*(M)$ by linearity over $\La$, and is associative. 
Moreover,
it  preserves the $\R$-grading.  

 This product $*_M$ gives $QH_{*}(M)$
the structure of a
graded commutative ring with unit $\1 = [M]$. Further, the invertible
elements in $QH_{\ev}(M)$ form a commutative group  
$QH_{\ev}(M,\La)^\times$
that acts on   $QH_*(M)$ by  quantum multiplication.

Now consider the fibration $P_\la\to S^2$ constructed from a loop $\la$
as in \S\ref{ss:geomint}.
As noted in~\cite{LMP}, the manifold $P_\la$ carries two canonical 
cohomology classes,
the first Chern class of the vertical
tangent bundle 
$$ 
c_{vert} = c_1(TP_\la^{vert})  \in H^2(P_{\la},\Z),
$$
and the   coupling class  $u_\la$, i.e.  the unique class
in $H^2(P_\la,\R)$ such that 
$$
i^*(u_\la) = [\om],\qquad u_\la^{n+1} = 0,  
$$
where $i: M\to P_\la$ is the obvious inclusion.

The next step is to choose a canonical (generalized) section class\footnote
 {
 By section class, we mean one that projects onto the positive generator
of $H_{2}(S^{2}, \Z)$.}  
 $\si_{\la}\in H_{2}(P_{\la}, 
 \R)/\!\sim$. In the general case, 
when $c_{1}$ and $[\om]$ induce linearly 
independent homomorphisms $H_2^S(M) \to \R$, $\si_\la$ is
defined by the requirement that
\begin{equation}\label{eq:seccl}
c_{vert}(\si_{\la}) \;= u_{\la}(\si_{\la})\; =\;0,
\end{equation}
which has a unique solution modulo the given equivalence.
We show in~\cite{Mcq} that when $M$ is weakly exact
such a class $\si_\la$ still exists and moreover is integral.
(The proof is included in Lemma~\ref{le:split} below.)
In the remaining  spherically monotone case, 
we choose $\si_\la$ so that $c_{vert}(\si_\la) = 0$. 

We then set
\begin{equation}\label{eq:qm2}
\Psi(\la) = \sum_{B\in \Hh} a_B\otimes e^{B}
\end{equation}
where, for all $c\in H_{*}(M)$,
\begin{equation}\label{eq:qm3}
a_{B}\cdot_{M} c = n_{P_{\la}}([M], [M], c\,; \si_{\la} - B).
\end{equation}
Note that $\Psi(\la)$ belongs to the strictly commutative part 
$QH_{ev}$ of $QH_{*}(M)$.
Moreover $\deg(\Psi(\la)) = 2n$ because $c_{vert}(\si_{\la}) = 0.$
It is shown in~\cite{Mcq} (using ideas from~\cite{Seid,LMP}) that
for all $[\la_1],[\la_2]\in \pi_1(\Ham(M))$
$$
\Psi(\la_1 + \la_2) = \Psi(\la_1)*\Psi(\la_2),\qquad 
\Psi(0) = \1,
$$
where $0$ denotes the constant loop.  Therefore 
$\Psi(\la)$ is invertible for all $\la$ and we get a 
representation
\footnote
{
$\Psi$ is 
called $\rho$ in~\cite{Mcq}, but we have changed the notation here, 
for obvious reasons.  Note also that the formula for $\Psi$ has been written
using sign conventions for $B$ that are different from those
in~\cite{LMP,Mcq}, to clarify
the inequalities  considered later.}
$$
 \Psi: \pi_{1}(\Ham(M,\om))\;\;\to\;\; (QH_{\ev}(M,\La))^{\times}.
$$
Since all $\om$-compatible forms are deformation 
equivalent, $\Psi$ is independent of the choice of $\Om$.

This is a mild reformulation of Seidel's representation 
from~\cite{Seid,Seid2}.  Our choice of canonical section class 
allows us to define $\Psi$ on the group
$ \pi_{1}(\Ham(M))$ itself rather than on the extension considered by Seidel.
The cost is that we have to allow $B$ to range in the real group $\Hh$
rather than in its integral lattice.   For 
 $\si_{\la}$ might have real coefficients,  and if it does, since we 
need to sum over classes $B$ such that $\si_{\la} - B$ is integral, 
we must allow $B$ to have real coefficients.   

Note also that the fact that $\Psi(\la) \ne 0$ implies that the 
fibration $P_{\la}\to S^2$ has a section (which, moreover, can be taken to be 
symplectic); equivalently, the loop $t\mapsto \la_t(x)$  contracts 
in $M$ for all $x\in M$.  This is an easy consequence of the proof of the Arnold 
conjecture: see for example~\cite{LMPf}\S1.3.   However, it can be 
proved more simply,  using only the compactness theorem and not gluing, by 
considering the limit of $J$-holomorphic sections in class $[pt \times 
S^2]$ of the trivial 
bundle  $M\times S^2 = P_\la \# P_{-\la}$ as this space degenerates into the 
singular union $P_{\la}\cup_M P_{-\la}$ as described in~\cite{Mcq}~\S2.3.2.

\subsection{Using $\Psi$ to estimate area}\label{ss:psiarea}

Consider the fibration $\pi: (P_{\la}, \Om)\to S^{2}$ corresponding to the 
class $[\la]\in \pi_{1}(\Ham(M))$.  Here $\Om$ is any $\om$-compatible
symplectic form on $P = P_{\la}$.  Hence, as remarked in \S\ref{ss:geomint},
 its cohomology class has the form
$$
[\Om] = u_{\la} + \pi^{*}([\al])
$$ 
where  $\area(P_{\la}, \Om)  =  \int_{S^2} \al$.

Following Seidel, consider the valuation $v:QH_{*}(M) \to \R$ defined by
\begin{equation}\label{eq:val}
v(\sum_{B\in \Hh} a_{B}\otimes e^{B}) = \sup \{\om(B): a_{B}\ne 
0\}).
\end{equation}
It follows from the definition of the quantum intersection product
in~(\ref{eq:qm0}), (\ref{eq:qm}) that
$v(a*b) \le v(a) + v(b).$  As Seidel points out, the following stronger statement is 
true.  Here we set
$$
\hbar = \hbar(M) = \min(\{\om(B)> 0: \mbox{ some }\;n_M(a,b,c;B)\ne 0\},
$$
 and  $\cap $ denotes the usual intersection product, so that $a*b - 
a\cap b$ is the quantum correction to the usual product.  Note that if all 
the invariants $n_M(a,b,c;B)$ with $a,b,c\in H_*(M)$ and 
$B\ne 0$ vanish, then $\hbar = \infty$.

\begin{lemma}\label{le:hbar}  For all $a,b\in
QH_*(M)$,  $
v(a*b - a\cap b) \le  v(a) +v(b) - \hbar(M)$. 
\end{lemma}
\proof{}  This follows immediately from the definitions.\QED

Seidel's results in~\cite{Seid2} are based on the following observation.

\begin{prop}[Seidel]\label{prop:main}  Suppose that $(M, \om)$ is not 
spherically monotone.  Then, for each loop $\la$ in $\Ham(M)$  
$$
 \area(P_{\la}, \Om) \;\;\ge\;\; v(\Psi(\la)).
$$
\end{prop}
\proof{}  Since $[\om]$ and $c_{1}$ are linearly independent on 
$H_2^{S}(M)$ we may define $\Psi$ using a section $\si_{\la}$ that 
satisfies~(\ref{eq:seccl}).
Let $\Psi(\la) = \sum_{B} a_{B}\otimes e^{B}.$
Then by definition $a_{B}$ is derived from 
a count of $J$-holomorphic curves in 
$(P_{\la}, \Om)$ in the class $\si_{\la} -B$.  If $a_{B}\ne 0$ then 
this moduli space cannot be empty. Hence
$$
\begin{array}{lcl}
0 \;\;< \;\;[\Om](\si_{\la} -B)  & = &  \pi^{*}([\al])(\si_{\la}) - \om
(B) \\
& = &  \int_{S^{2}}\al - \om(B) \\
&=& \area(P_{\la}, \Om) - \om(B),
\end{array}
$$
as required.\QED

\begin{cor}\label{cor:main} In these circumstances $\ell^+(\la) \ge v(\Psi(\la)).$
\end{cor}
\proof{}  Since $\ell^+(\la) = \trho^+([\la])$, this follows by combining
Lemma~\ref{le:polt2} with Proposition~\ref{prop:main}.\QED
 
We will apply this in~\S\ref{sec:exam} to calculate the lengths of loops.  
Note that the above argument applies when $\Psi$ is normalized 
using $\si$ for which $u_{\la}(\si) = 0$: the value of 
$c_{vert}(\si)$ is irrelevant here.  Hence it can apply
in the spherically monotone case when $[\om] = \ka c_{1}$ on 
 $\pi_{2}(M)$ for $\ka \ne 0$.

 \begin{rmk}\label{rmk:pol}\rm
Another closely related way of getting an estimate 
in the monotone case when $[\om] = \ka c_{1}$ for $\ka > 0$ was 
observed by  Polterovich in~\cite{Pac} Thm 2.A.
He defined a 
homomorphism $I:\pi_{1}(\Ham(M, \om))\to \R$ as follows: for each 
$\la$
choose a section class $\si'_{\la}$ in $P_{\la}$ 
such that $c_{vert}(\si'_{\la}) = 0$ and then set
$$
I(\la) = u_{\la}(\si'_{\la}).
$$
A dimension count shows that the only section classes $\si$ that contribute
to $\Psi(\la)$ have $c_{{vert}}(\si)\le 0$.  Hence, 
given any such $\si$ and  any symplectic 
form $\Om$ on $P_{\la}$, 
$$
0 \;\;<\;\; \Om(\si) \;\;=\;\; 
\area(P_{\la}, \Om) + u_{\la}(\si) \;\;= \;\;\area(P_{\la}, \Om) + 
u_{\la}(\si'_{\la} + B)\;\; \le \;\; \area(P_{\la}, \Om) + I(\la),
$$
since $u_{\la}(B) = \ka c_{1}(B)$ for $\ka > 0$ is nonpositive.
Therefore, since $I(-\la) = -I(\la)$,
$$
\trho_{f}(\la)\;\; =\;\; \inf(\area(P_{\la}, \Om)) + \inf(\area(P_{-\la}, 
\Om)) \;\;\ge \;\;|I(\la)|.
$$
This estimate is weaker than the previous one if $c_{1}(B) < 0$.  
However, because
$I$ is a homomorphism, one immediately finds $\trho_{f}(k\la)\ge k|I(\la)|$.
See also~\cite{Seid2}.
\end{rmk}

\subsection{The nonsqueezing theorem}\label{ss:ns}

Let $QH_+(M)$ denote the set 
$$
QH_+(M) = \{x \in \sum_{i\le 2n-2} H_{i}(M)\otimes \La\subset QH_{*}(M)\}.
$$
We first give a simple criterion for the nonsqueezing theorem to hold.
Recall from above that every fibration $P_{\la}\to S^{2}$ admits a 
generalized section class $\si_{\la}$ on which the coupling class
$u_{\la}$ vanishes, except possibly when $(M, \om)$ is weakly exact
and $c_{1}$ does not vanish on $H_{2}^{S}(M)$.

\begin{lemma}\label{le:ns}  Suppose that $(P_{\la}, \Om)\to S^{2}$ admits a 
generalized section class $\si_{\la}$ on which 
$u_{\la}$  vanishes, and that the corresponding element $\Psi(\la)$ has the 
form
$$
\Psi(\la) = \1\otimes \mu + x,
$$
where $x \in QH_+(M)$ and $\mu = \sum q_{B}\,e^{B}$ has some nonzero 
coefficient $q_{B}$ with $\om(B) \ge 0$.
Then the nonsqueezing theorem holds for $(P_{\la}, \Om)$.
\end{lemma}
\proof{}  The hypotheses imply that 
 $ n_{P}([M], [M], pt; \si_{\la} - B) = q_{B}\ne 0$. 
Since this invariant counts perturbed $J$-holomorphic stable maps in class 
$\si_{\la} - B$ through an arbitrary point, it follows that there is such 
a curve through every point in $P$.  Since the perturbation can be taken 
arbitrarily small, it follows from Gromov's compactness theorem that there 
has to be some $J$-holomorphic stable map in this class through every point 
in $P$.
 Hence the usual arguments (cf~\cite{LMe}, for example)
 imply that the radius  $r$ of any embedded ball satisfies the 
 inequality:
 $$
 \pi r^{2} \;\;\le\;\; [\Om](\si_{\la} - B) \;\; \le 
\;\;[\Om](\si_{\la})\;\; =\;\; \area(P_{\la}, \Om).
 $$
The result follows.\QED

It will be convenient to make the following definition.

\begin{defn}\label{def:sns} We  say that the fibration $(P, 
\Om)\to S^{2}$  with fiber $M$ has  
 a {\bf good section of weight} $\ka$ if there is 
a class $\si_{A} \in H_{2}(P)$ such that
\begin{itemize}
\item[(i)] $ n_{P}([M], [M], pt; \si_{A}) \ne 0$;
\item[(ii)]  $u(\si_{A}) = - \ka$ where $u$ is the coupling class.
\end{itemize}
Note that $\ka$ could be positive or negative.
\end{defn}

The previous lemma shows that any fibration with a good section
of weight $0$
has the nonsqueezing property.  More generally, the same argument proves 
the following weighted nonsqueezing property.

\begin{lemma}\label{le:ns1} 
Suppose that $(P_{\la},\Om)$ has a good section of weight $\ka$.
  Then the radius $r$ of an 
embedded ball in $(P,\Om)$ is constrained by the inequality:
$$
\pi r^2 \le \area\,(P, \Om) - \ka.
$$
In particular, if $\Psi$ is defined relative to a section class $\si_{\la}$
on which $u_\la$ vanishes, we may take $\ka$ to be the maximum
of $\om(B)$ where $q_B\ne 0$ in the expression for $\Psi(\la)$.
\end{lemma}

It is not hard to find conditions under which $
\Psi(\la) = \1\otimes \mu + x,
$ where $\mu \ne 0$.  The tricky point is to 
find ways of estimating the maximal weight of a good section.
In this subsection we will describe situations in which there is a good 
section of weight $0$ so that
the nonsqueezing theorem holds, leaving the discussion of the more general 
case to~\S\ref{ss:nsgeo}.

To find good sections of weight $0$, we can
use arguments  in~\cite{Mcq} that give conditions 
under which the (usual) cohomology 
ring $H^{*}(P)$ splits.  The idea is the following.  

Suppose that for some $A\in H_{2}(M)$ 
we can define a map $s_{A}: H_{*}(M) \to H_{*+2}(P)$
such that
\begin{equation}\label{eq:sA}
s_{A}(pt) = \si_{\la} - A, \quad s_{A}(a)\cap [M] = a,\quad 
s_{A}(a\cap b) = s_{A}(a)\cap s_{A}(b),
\end{equation}
for all $a,b\in H_{*}(M)$.
Then the Poincar\'e dual map
$$
r: H^{*}(M) \to H^{*}(P), \quad \al\mapsto \PD_{P}(s_{A}(\PD_{M}(\al)))
$$
is a ring homomorphism such that $\io^{*}\circ r =\id$, where 
$\io:M\to P$ is the inclusion.
In particular,  $r([\om]) = u_{\la}$, since 
$r([\om])^{n+1} = 0$ and $u_{\la}$ is the unique extension of $\om$ 
such that $u_{\la}^{n+1} = 0$.  Further,
$$
\PD_{P}(u_{\la}^{n}) = s_{A}(\PD_{M}[\om]^{n}) = s_{A}(pt)/n! = 
(\si_{\la} - A)/n!.
$$
Hence $u_{\la}(\si_{\la} - A) = 0$.  Since $u_{\la}(\si_{\la}) = 0$ 
by construction, we have $\om(A) = 0$ as required.  Thus

\begin{lemma}\label{le:goods}  Suppose that there is a 
splitting $s_{A}: H_{*}(M) \to H_{*+2}(P)$ satisfying~(\ref{eq:sA}).
Then $\om(A) = 0$.
\end{lemma}

There are two ways to construct such a splitting $s_{A}$.  First suppose 
that there is an $\Om$-tame almost complex structure
$J$ on $(P, \Om)$ such that the moduli space 
$\Mm_{J}$ of  $J$-holomorphic 
curves of class $\si_{A} = \si_{\la} - A$  is  compact  
(where we assume them to be parametrized as sections) and of 
dimension $2n$.
Then,  there are evaluation 
maps
$$
e:\Mm_{J}\times S^2\to P,\qquad e_{0}:\Mm_{J}\to M
$$
of equal degree $q$.  If $q \ne 0$ we define
\begin{equation}\label{eq:goods}
s_{A}: H_{*}(M) \to H_{*+2}(P):\quad
a\;\;\mapsto\;\;\frac 1q\;e_{*} (e_{0}^{!}(a)\times 
[S^{2}]), 
\end{equation}
where  $(e_{0})^{!}: H_{*}(M)\to H_{*}(\Mm)$  is the homology 
transfer (defined as the Poincar\'e dual $\PD_{\Mm} \circ e_{0}^{*} 
\circ\PD_{M}$
of the pullback in cohomology.)
More generally, if all we know is that the invariant
$$
n_{P}([M], [M], pt; \si_{A}) = q_{A}\ne 0,
$$
we define $s_{A}(a)$ to be the unique class in $H_{*}(P)$ such that
$$
s_{A}(a)\cdot_{P} v =\frac 1{q_A}  n_{P}(a, [M], v; \si_{A}),\qquad v\in H_{*}(P).
$$
(Here, as elsewhere, $H_{*}$ denotes rational homology.)  

The next lemma describes situations in which $s_A$ satisfies the conditions 
in~(\ref{eq:sA}). 
Statement (i) is due to Seidel\footnote
{Private communication.} and suffices to prove all our main results,
including Proposition~\ref{prop:ns} and Theorem~\ref{thm:geod}.
We include the proof because it is completely 
elementary, even though one could equally well argue using the other 
parts of the lemma below.
Note that (ii) is a corrected version 
of~\cite{Mcq} Proposition~3.4(i).\footnote
{ 
In the preprint version of~\cite{Mcq}, Proposition 3.4 is numbered as 3.21.}
  As in~\cite{Mcq}, the letters $a,b,c$ 
denote either elements of $H_{*}(M)$ or their images in $H_{*}(P)$,
and $u,v,w$ denote general elements in $H_{*}(P)$.  
Also $B\in H_{2}(M, \Z)/\!\sim\;= \Hh$.

\begin{lemma}\label{le:split}  With notation as above,
\begin{equation}\label{eq:cond}
s_{A}(a)\cap [M] = a,\quad s_{A}(a\cap b) = s_{A}(a)\cap s_{A}(b), 
\;\;\mbox{for all }\; a,b\in H_{*}(M)
\end{equation}
under each of the following conditions:
\begin{itemize}
\item[(i)]  $s_{A}$ is defined as above from a compact moduli
space  $\Mm: = \Mm_{J}$.
\item[(ii)]  The only nonzero Gromov--Witten invariants of the form
$
n_{P}([M], a, v; \si_{A} - B),$ with  
$\om(B) \ge 0$   have $B = 0$.
\item[(iii)]  The $3$-point invariants $n_{M}(a,b,c; B), B\ne 0,$ 
vanish and $n_{P}([M], [M], pt; \si_{A} - B)= 0$  when 
$\om(B)> 0$.
\item[(iv)]
All  $3$-point vertical invariants $n_{P}(u,v,w; B)$, $B\in H_{2}(M) - 
\{0\}$, in 
$P$ vanish,  as do all  
$4$-point invariants $n_{M}(a,b,c,d; B)$, $B\ne 0$, in $M$.
\end{itemize}
\end{lemma}
\proof{} (i)  It is slightly easier to prove the cohomological version of 
(i).  However, we phrase the argument in homology to make this argument 
closer to the proof of (ii) -- (iv).

If $e:X\to Y$ is any map between two closed
manifolds of the same 
dimension, then  the homology transfer 
$$
e^!: H_*(Y)\to H_{*}(X):\quad a\mapsto \PD_X(e^*(\PD_Y a)) 
$$
 has the following properties: 
$$
\begin{array}{lclll}
(a) &\quad& e^!(v\cap w) & = &  e^!(v)\cap e^!(w),\\
(b) & \quad& e_*(v)\cap w & = & e_*(v \cap e^! (w)), \\
(c) &\quad & e_*(v) \cdot  w & = &  v\cdot e^!(w),\\
(d) & \quad & e_* e^!(v)  & = &  ({ deg\,} e)\, v,
\end{array}
$$
Further, if  if $e_0:\Mm \to M, e:\Mm\times S^2\to P$
 are as above, it is not hard to
check that
$$
(e)\qquad\qquad  e^!(v)\cap [\Mm] =  e_0^{\,!} (v\cap [M]).
$$
Thus, by~(\ref{eq:goods})
\begin{eqnarray*}
s_{A}(a)\cap [M]   & = & \frac 1 {q}\, e_*\left(e_0^{\,!}(a)\times 
[S^2]\right)
\cap [M]\\
& = & \frac 1 {q} \,e_*\left(e_0^{\,!}(a)\times [pt]\right)\\
& = & \frac 1 {q}\, (e_0)_* e_0^{\,!}(a) \, =\, a.
\end{eqnarray*}
It remains to prove the second half of
condition~(\ref{eq:cond}).  In view of (a) above, this would be 
obvious if
$e_*$ respected the cap product.  Since this is not the case, we must 
use a different approach.

Let us write $e^!(s_{A}(a)) = [a_0\times S^2] + \io(a_1)$ where $a_1 \in
H_*(\Mm)$ and $\io:H_*(\Mm)\to H_*(\Mm\times S^2)$ is induced by the 
inclusion.  It follows easily from (e) above that $a_0 =
e_0^{\,!}(a).$  Hence, by (d), $(e_0)_*(a_1) = 0$.  By (c), 
this means that
$a_1\cdot e_0^{\,!}(b) = 0$ for all $b$.  Hence 
\begin{eqnarray*}
q\,s_{A}(b)\cdot s_{A}(a)  & = &  e_*([e_0^{\,!}(b)\times S^2])\cdot s_{A}(a)\\
& = & [e_0^{\,!}(b)\times S^2]\cdot (a_0\times [S^2] + \io(a_1))\\
& = & a_1\cdot e_0^{\,!}(b) \, = \,0
\end{eqnarray*}
for all $a,b$.  Again using  $s_{A}(a)\cap [M] = a$, we have that
$$
s_{A}(a)\cap s_{A}(b) = s_{A}(a\cap b) + \io(x)
$$
for some $x\in H_*(M)$.  Further $x = 0$ if and only if 
$$
(s_{A}(a)\cap s_{A}(b))\cdot_P s_{A}(c) = \io(x)\cdot_P s_{A}(c) = x\cdot_M c = 0
$$
for all $c$.  But, as above,
\begin{eqnarray*}
& &q\, s_{A}(c)\cdot_P (s_{A}(a)\cap s_{A}(b))\\
& & \qquad\qquad =  [e_0^{\,!}(c)\times S^2] \cdot 
\left(e^!(s_{A}(a))\cap e^!(s_{A}(b))\right)\\
& & \qquad\qquad =   [e_0^{\,!}(c)\times S^2] \cdot 
\left([e_0^{\,!}(a)\times S^2]
\cap \io(b_1) +  \io(a_1)\cap [e_0^{\,!}(b)\times S^2]\right) \\
& & \qquad\qquad =  e_0^{\,!}(c\cap a)\cdot b_1  \;\;\pm\;\;  e_0^{\,!}(c\cap 
b)\cdot a_1\\
& & \qquad\qquad =  0.
\end{eqnarray*}
This completes the proof of (i).\QED

The proofs of (ii), (iii) and (iv) can be found in~\cite{Mcq};
(iv) is Proposition~3.4, while (ii) is a corrected form of
Proposition~3.5~(i).   We have added the extra assumption that there 
are no nonzero invariants of the form
$
n_{P}(v, [M], a; \si_{A} - B)$ with $B\ne 0$ in order for the proofs 
of \cite{Mcq}~Lemmas 3.10 and 3.11 to hold.  (iii) is a slightly more 
general version
of ~\cite{Mcq}~Proposition~3.5~(ii). We now allow invariants of the
form $n_P([M], [M], pt; \si_{A} - B)$ with $\om(B) < 0$ to be 
nonzero.  But clearly these do not contribute to the sums considered
in ~\cite{Mcq}~Lemmas 3.9, 3.10, 3.11.
\QED

\begin{cor}\label{cor:split}  
Consider a  fibration $(P_{\la}, \Om) \to S^{2}$
 such that
$n_{P}([M], [M], pt; \si_{A}) \ne 0$ for
some class $\si_{A}\in H_{2}(P)$.  If  one of the  conditions
in Lemma~\ref{le:split} also holds
then $\si_A$ is a good section of 
weight $0$ in $(P_{\la}, \Om)$.
\end{cor}
\proof{} This holds by  Lemma~\ref{le:goods}.
Observe that the conclusion holds even when we cannot assume that 
$u_{\la}(\si_{\la}) = 0$ since we prove in all cases that  
$u_{\la}(\si_A) = 0$. \QED

With these preliminaries, we are now ready to prove
 Proposition~\ref{prop:ns}.  In view of Lemma~\ref{le:ns}
 it is an immediate consequence of the 
next result. Here 
$$
\La^-\;\; =\;\; \{\mu\in \La\;:\; \mu = q_{0} e^{0} + \sum q_{C} e^{C}, \;\; \om(C) < 0, \; 
q_{C}\in \Q, \;q_{0}\ne 0.\}
$$

\begin{prop}\label{prop:ns2}  Suppose that $(M, \om)$ is spherically 
rational with index of rationality $q(M)$ and
let $\la$ be a loop in $\Ham(M,\om)$.
If
$\ell^{+}(\la) + \ell^{-}(\la) < q(M)$ then $\Psi(\la)$ and 
$\Psi(-\la)$ both have the form $\1\otimes \mu + x$ with $\mu\in 
\La^{-}$,
$x \in QH_{+}(M)$.  In particular, they both have a good section of 
weight $0$.
\end{prop}
\proof{}  We will assume that $M$ is not weakly exact since this case is 
dealt with in Lemma~\ref{le:Jsph}.  Hence we can choose $\si_\la$ so that 
$u_\la(\si_\la) = 0$.

 Choose $\eps > 0$ so that
$\ell^{+}(\la) + \ell^{-}(\la) < q(M) - 2\eps.$
By  Proposition~\ref{prop:polt}, there is a $\om$-compatible symplectic form
$\Om_{\la}$ on $P_{\la}$ with area $< \ell^{+}(\la) + \eps$, and  a 
similar form $\Om_{-\la}$ on $P_{-\la}$ with area $< \ell^{-}(\la) + \eps$.

Write
$$
\Psi(\la) =  \sum_{B\in\Hh} \;\1\otimes\,q_{B} e^{B} + x,
\quad
\Psi(-\la) =  \sum_{B'\in\Hh}\; \1\otimes \,q_{B'}'\,e^{B'} + x'
$$
where $q_{B},q_{B'}'\in \Q$ and $x,x'\in QH^{+}$.
Note first that for all  $e^{B}$ (resp. $e^{B'}$) that occur in $\Psi(\la)$ 
(resp.  $\Psi(-\la)$) with
nonzero coefficient,
$$
\om(B) < \ell^{+}(\la) + \eps,\qquad \om(B') < \ell^{-}(\la) + \eps.
$$
This  holds from the definition of the 
coefficients via Gromov--Witten 
invariants $n_{P}([M], [M], c; A)$ where $A =\si_{\la} - B$ 
(resp. $\si_{-\la} - B$), and the fact that 
$\Om_{\la}$ (resp. $\Om_{-\la}$) must have positive integral on $A$:
see equations~(\ref{eq:qm2}), (\ref{eq:qm3}).

Next apply the valuation $v$ in~(\ref{eq:val}) to the identity
$$
\Psi(\la)*\Psi(-\la) = \Psi(0) = \1.
$$
We claim that at least one of $\Psi(\la),  \Psi(-\la)$
has a term $\1\otimes q_{B} e^{B}$
 with $q_{B}\ne 0, \om(B) \ge 0$.
For otherwise the product $x*x'$ must contain the term  $\1\otimes 
e^{0}$ with a nonzero 
coefficient.  Because this term appears in $x*x' - x\cap x'$ we find
from Lemma~\ref{le:hbar} that
$$
0 \;=\; v(\1\otimes e^{0}) \;\;\le\;\; v(\Psi(\la)) + v(\Psi(-\la)) - q(M)\;\; \le\;\; 
\ell^{+}(\la) + \ell^{-}(\la) + 2\eps  - q(M) \;\;<\;\; 0,
$$
a contradiction.

 Therefore, replacing $\la$ by $-\la$ if necessary we may suppose that
 $$
 \Psi(\la) =  \1\otimes \mu e^{A} + x 
 $$
where $x\in QH_{+}(M), 
 0 \le \om(A) <q(M) - \eps,$ and $ \mu\in \La^{-}$.
 The lemma will follow if we show that  $\om(A) = 0$.

 To do this, consider the class $\si_{A} = \si_{\la} - A$ as above.  
 In view of Corollary~\ref{cor:split} it suffices to prove the 
 following claim. 
\MS

\NI
{\bf Claim:}\,\, {\it There is an $\Om$-tame $J$ on $P$
such that the moduli space $\Mm_{J}$ of unparametrized 
$J$-holomorphic 
curves in $(P, \Om)$ of class $\si_{A}$ is a compact
manifold of dimension $2n$.}\MS

\NI
{\bf Proof of Claim:}\,\,  We first show that $\Mm_J$ is compact for all 
fibered $J$. (Recall that an almost 
 complex structure $J$ on a symplectically 
 fibered space $P\to S^{2}$ is called 
 fibered if each fiber is $J$-holomorphic.)
If not, there is a sequence of $\si_{A}$-curves that converges to a 
stable map.  One component of this stable 
map will be a section and the others will each lie entirely in a fiber.
There must be at least one bubble in a fiber, which will use up a 
minimum of $q(M)$ in energy.   Since 
$$
[\Om](\si_A) = \eps +  \om(A) < q(M)
$$
 this is impossible.  Next observe that the curves in $\Mm_J$ are all 
 embedded so that they can be regularized by choosing a generic $J$. 
It is not hard to see that this can be chosen to be fibered: compare 
Lemma~4.3 of~\cite{Mcq}.\footnote
{
Lemma 4.9 in the preprint.}   Alternatively, use the compactness theorem again 
to conclude that $\Mm_J$ is compact for
every $J$ that is sufficiently close to a fibered $J$, and then choose a 
regular one from among these.  This proves the claim and hence the 
Proposition.\QED

\begin{rmk}\rm Using~\cite{Mcq} Proposition 3.4 one can conclude from the 
above argument
 that   if $(M, \om)$ is spherically rational then
 whenever $[\la]\in \pi_{1}(\Ham(M))$ and its inverse $[-\la]$ 
 have sufficiently 
$\trho$-short representatives  they are in the kernel of 
the homomorphisms $\oI_{c}$ and $\oI_{u}$ 
in~\cite{LMP,Mcq}.  Moreover the rational cohomology rings $H^{*}(P_{\la})$ 
and   $H^{*}(P_{-\la})$ split as products.
\end{rmk}  

Further extensions of the proof of Proposition~\ref{prop:ns2} are 
discussed in \S\ref{ss:nsgeo}.  We end this section by
 proving Proposition~\ref{prop:nsa}. 

\begin{lemma}\label{le:Jsph}\begin{itemize}
\item[(i)]  When $c_1= 0$ on $\pi_2(M)$ and all
 $3$-point Gromov--Witten invariants on $M$ vanish, 
 every fibration $(P, \Om) \to S^{2}$ has a good section of weight $0$. 
\item[(ii)]  The same statement holds if $(M, \om)$ is weakly exact.
\item[(iii)] If $M = \CP^n$
the $(n+1)$st multiple $(n+1)[\la]$ of each loop $\la$
has a good section of weight $0$.
\end{itemize}
\end{lemma}
\proof{} (i) Assume first that $c_1$ vanishes on $\pi_2(M)$ 
but $\om$ does not, and choose $\si_{\la}'$ so that
$u_{\la}(\si_{\la}') = 0$. We claim that $c_{vert}(\si_{\la}')$ must
also be zero.  For
$c_{vert}$ takes the same value on all section classes $\si_{\la}' - 
B$,
and since at least one invariant $n_{P}([M], [M], a;\si_{\la}' -B)\ne 0$
this value must be $\le 0$.  On the other hand
$$
c_{vert}(\si_{\la}') + c_{vert}(\si_{-\la}') \;\; = \;\; 
c_{vert}(\si_{\la}'\#\si_{-\la}') \;\; = \;\;  0,
$$
since the concatenation $\si_{\la}'\#\si_{-\la}'$ is the canonical 
section in the product $M\times S^{2}$: see~\cite{LMP}.
Hence $c_{vert}(\si_{\la}') = 0$ so that $\si_{\la}' = \si_{\la}$,
and $\Psi(\la) = \1\otimes \mu$ for 
some $\mu \in \La$. Therefore, there are classes $B$ so that
$n_{P}([M], [M], pt; \si_{\la} - B) \ne 0$ and we choose $A$ from 
among them so that $\om(A) $ is a maximum.
This means that condition (iii) in Lemma~\ref{le:split} holds.  
Therefore the result follows from Corollary~\ref{cor:split}.\QED

To prove (ii), suppose that $(M, \om)$ is weakly exact.  We said before that by
the results of~\cite{Mcq} 
 we may choose $\si_{\la}$ so that both classes $c_{vert}$ 
and $u_{\la}$ vanish on it.  For the sake of completeness, we will 
give this argument now.  Suppose first that  $c_{1}$ vanishes on 
$H_{2}^{S}(M)$ as well.  Then there is only one section class (up to 
equivalence) and we call it $\si_{\la}$. 
The argument in (i) above shows that $c_{vert}$ must 
vanish on this class.  Further, the moduli space $\Mm_{J}$ of curves 
in this class must be compact (since there are no $J$-holomorphic 
curves in $M$).  Hence $u_{\la}(\si_{\la}) = 0$ as well.  Since 
$\Psi(\la) = q\1 \ne 0$ the result follows.

When $c_{1} \ne 0$ on $H_{2}^{S}(M)$, we choose $\si_{\la}$ so that
$c_{vert}(\si_{\la}) = 0$.  As in (i), we choose $A$ so that $\om(A)$ 
is maximal among the classes $B$ for which
$n_{P}([M], [M], pt; \si_{\la} - B) \ne 0$.  The argument then 
continues either as in (i) or as in the proof of 
Proposition~\ref{prop:ns2}.  In particular, it shows that $u_{\la}(\si_{A}) = 
0$.\QED

When  $M = \CP^n$, we again choose $\si_\la$ so that $c_{vert}(\si_\la) = 0$,
so there is no control over $u_{\la}(\si_\la)$.
Because $c_1(L)= n+1$, where $L= [\CP^1]$,
and  $-n \le c_{vert}(\si_\la -B) \le 0$ whenever
\begin{equation}\label{eq:gw}
n_P([M], [M], c; \si_\la - B)\ne 0,\qquad c\in H_{*}(M),
\end{equation}
there can be at most  one
class $B = r L$ with nonzero invariant.  
Hence $\Psi(\la)$ has the form $a \otimes 
e^{r L}$, where  $a \in H_{2k}(\CP^n)$ is homogeneous.  
Since $\deg \Psi(\la)= 2n$ we must have $r(n+1) = n-k$.
Therefore 
$$
\Psi((n+1)\la) = a^{(n+1)}\otimes e^{(n-k) L} = q \1,\quad q\in \Q - 
\{0\},
$$
since  the hyperplane class  $b\in QH_{2n-2}(\CP^{n})$ satisfies the
relation $b^{(n+1)} = \1\otimes e^{-L}$ and $a = q' b^{n-k}$.   
Therefore, if $Q\to S^{2}$ denotes the fibration corresponding to the 
loop $(n+1)\la$, 
$$
n_{Q}([M], [M], pt; \si_{(n+1)\la}) \ne 0.
$$
Again using the fact that $c_{1}(L) = n+1$ and taking $\si_{A} = 
\si_{(n+1)\la}$, we find that
condition (ii) in Lemma~\ref{le:split} holds.  Hence
$u_{(n+1)\la}(\si_{(n+1)\la}) = 0$ and the result follows.
Observe that in this case condition (i) in
Lemma~\ref{le:split} also holds for generic $J$ on $P$
since for each $k\ge 1$ the moduli space of curves in 
$P$ of class $\si_A - kL$ must vanish 
for reasons of dimension. \QED

\begin{rmk}\label{rmk:nsthm} \rm 
 \begin{itemize}
\item[(i)]  The above argument about $\CP^{n}$ 
 applies to any  manifold $M$ for which $\Psi(\pi_{1}(\Ham))$ is a finite
 subgroup of  $QH_{ev}(M)^{\times}$.  Here we should assume that $M$ 
 is not spherically monotone so that we always have $u_{\la}(\si_{\la}) = 0$.  
 This was unnecessary for  $M = \CP^{n}$ since
 $c_{1}(L)$ is so large.
 \item[(ii)]  Lemma~\ref{le:ns1} shows that
there is some nonsqueezing inequality for any fibration 
$(P_{\la},\Om)\to S^{2}$ that has good sections, i.e. has
$\Psi(\la) = \1\otimes \mu + x$ 
for some $\mu \ne 0$.  One cannot say anything for general fibrations
$P_{\la}\to S^{2}$  unless the embedded ball 
is disjoint from one fiber.  In the latter case, one
can make $P_{\la}$ symplectically trivial by taking the 
fiber sum with $P_{-\la}$, and deduce that 
the radius $r$ of any symplectic ball in $(P_{\la}, \Om)$ 
that misses a fiber
is constrained by the inequality
$$
\pi r^{2}\le \area(P_{\la}, \Om) + \ell^{-}(\la).
$$
More generally, as  suggested by Polterovich,\footnote
{Private communication.}
one can  consider nonsqueezing for fibrations 
$$
(P_{\tphi}, \Om)\to D
$$
 with nontrivial  boundary monodromy $\tphi$.  
Completing these to symplectically trivial fibrations 
over $S^{2}$ by adding a fibration with boundary monodromy 
$\tphi^{-1}$, one finds that the radius $r$ of 
any symplectic ball in $(P_{\tphi}, \Om)$
satisfies 
$$
\pi r^{2} \le \area(P_{\tphi}, \Om) + \trho^{-}(\tphi).
$$
Further, if the nonsqueezing theorem holds for all loops in 
$\Ham(M)$, one has the inequality
$$
\pi r^{2} \le \area(P_{\phi}, \Om) + \rho^{-}(\phi).
$$
\end{itemize}
\end{rmk}

\subsection{Weighted nonsqueezing and geodesics}\label{ss:nsgeo}

In this section we  prove Proposition~\ref{prop:ns3} and 
Corollary~\ref{cor:nsgeo}, and hence
complete the proof of Theorem~\ref{thm:geod}.\MS

\NI
{\bf Proof of Proposition~\ref{prop:ns3}}.\,\,
By hypothesis there are fibrations $(P_{\pm\la}, \Om)$ with area $< 
\hbar/2$. Apply the valuation $v$ in~(\ref{eq:val}) to the identity
$$
\Psi(\la)*\Psi(-\la) = \Psi(0) = \1,
$$
as in the proof of Proposition~\ref{prop:ns2}.  As there, because 
$\area(P_{\pm\la}, \Om)< \hbar/2$, the product
$x*x'$ cannot contain the term  $\1\otimes e^{0}$ with a nonzero 
coefficient.  Therefore for $\la'$ equal to at least one of $\la$ or $-\la$, 
$\Psi(\la')$
has a term $q_{B}\1\otimes  e^{B}$
 with $q_{B}\ne 0$ and $0 \le  \om(B) < \area(P_{\la'}, \Om)$.  

Set 
$$
\eps(\la') = \max \{\om(B): q_{B}\ne 0 \;\mbox{ in }\;\Psi(\la')\},\quad
\la' = \pm\la.
$$
 The equation $\Psi(\la)*\Psi(-\la) = \1$ implies that
$\eps(\la) = -\eps(-\la)$.  Moreover,  
 by Lemma~\ref{le:ns1}, the radius $r$ of any embedded ball in 
 $(P_{\la'}, \Om)$ 
satisfies
$$
\pi r^2 \le \area(P_{\la'}, \Om) - \eps(\la').
$$
Hence we may take $\eps = -\eps(\la)$.\QED

\begin{rmk}\rm
The above argument uses only the first half of the proof of 
Proposition~\ref{prop:ns2} since it is not clear what  hypothesis 
would guarantee that the Claim holds, i.e. that an appropriate  moduli 
space of sections is compact.   Although for each individual 
$\om$-tame $J$ the minimum energy $\hbar(M,J)$ of a nontrivial 
$J$-holomorphic bubble is positive,  the 
minimum of $\hbar(M,J)$ over all such $J$ will  not in general 
be strictly positive since any symplectically embedded $2$-sphere is
$J$-holomorphic for some tame $J$.
We cannot restrict ourselves to a compact set of $J$ since we must 
consider all loops $\la$ in $\pi_{1}(\Ham)$, each of which gives rise 
to some  $2$-parameter family of $J$ on the fibers of $P_{\la}$.  
One might also try  to establish nonsqueezing by 
using  Lemma~\ref{le:split} $(iv)$.
For this one would 
  need a constant $\hbar>0$ such that all vertical Gromov--Witten invariants
$n_{P}(u,v,w; B)$ vanish,  where  $u,v,w$ are arbitrary elements in 
$H_{*}(P_{\la})$ and 
 $0 <\om(B) < \hbar$.  Again, because we are 
considering the Gromov--Witten invariants of an arbitrary fibration 
$P_\la$ rather than those of a compact manifold, it is not 
clear that such $\hbar$ exists.  Therefore, our present methods do not suffice 
to show that 
the nonsqueezing theorem holds for all $\trho$-short loops.
\end{rmk}

\NI
{\bf Proof of Corollary~\ref{cor:nsgeo}}.\,\,
If $(\tnu, \nu)$ is one of the pairs
$(\trho,\rho), (\trho_f,\rho_f)$,  we must show
that if a path $\ga$ has $\tnu$-length $< \hbar(M)/4$, 
minimizes $\tnu$ in $\THam(M)$ and 
is sufficiently $C^2$ close to the
identity, then
it is 
$\nu$-minimizing in $\Ham(M)$.  

As in the proof of Lemma~\ref{le:nsgeo} we will carry out the argument 
for the pair $(\trho_f,\rho_f)$.  
Suppose that $\ga$ does not minimize $\rho_{f}$.
Then there are paths 
 $\psi_{t}^{\pm}$ from $id$ to $\phi = \phi_1^H$  
generated by $K_{t}^{\pm}$ and $\de > 0$ such that
 $$
\Ll^{+}(K_{t}^{+}) + \Ll^{-}(K_{t}^{-}) \;= \;\Ll(\ga) - \de \;<\; \Ll(\ga) \;=\; 
\Ll^{+}(H_{t}) + \Ll^{-}(H_{t}).
$$ 
As before, we may assume that:
$$
\Ll^{+}(K_{t}^{+}) \;\; = \;\;  \Ll^{+}(H_{t}) - \de' 
\;\; < \;\; \Ll^{+}(H_{t}),\qquad
\Ll^{-}(K_{t}^{-}) \;\; = \;\; \Ll^{-}(H_{t}) - \de + \de',
$$
for some $\de' > 0$.   Let   $\la = \be * (-\ga)$ as before
so that $P_{K^{+},H} = P_\la$, $P_{H,K^{-}} = P_{-\la}$.
Then for small $\eps $
$$
\begin{array}{lcl}
\area(P_{\la}(\eps), \Om_{0}) & = &  \Ll(\ga)  - \de' + \eps \;\;<\;\; \Ll(\ga)
\;\; < \;\;
\hbar/4,\\ 
\area(P_{-\la}(\eps), \Om_0) & = & \Ll(\ga) - \de + \de' + \eps \;\;\le\;\;  
2\Ll(\ga) < \hbar/2.
\end{array}
$$
By Proposition~\ref{prop:ns3} there is $\ka$ with $|\ka| \le \hbar/2$
 such that embedded balls satisfy
$$
\pi r^2 \le \area(P_{\la}(\eps), \Om) + \ka,\quad 
\pi r^2 \le \area(P_{-\la}(\eps), \Om) - \ka.
$$
But, 
because $\ga$ is $C^2$-close to the
identity and must have fixed extrema, it follows from Lemma~\ref{le:1} that
both $(P_\la(\eps), \Om)$ and $(P_{-\la}(\eps), \Om)$ 
contain embedded balls of 
capacity $\pi r^2 = \Ll(\ga) > \area(P_\la(\eps), \Om)$.  Hence 
 $\ka > 0$.  Further,
$$
\begin{array}{lclcl}
\Ll(\ga) &\le & \area(P_\la(\eps), \Om) + \ka & = &\Ll(\ga) -\de' + \eps + 
\ka\\
\Ll(\ga) & \le & \area(P_{-\la}(\eps), \Om) - \ka & = &
\Ll(\ga) - \de + \de' + \eps - \ka.
\end{array}
$$
Adding, we find $0 \le -\de + 2\eps$.  Since $\de$ is positive and 
 $\eps$ can be arbitrarily small, this 
is impossible.  Hence result.\QED


\section{Estimating lengths of loops}\label{sec:exam}

Let $M_{*}$ be the one point blow-up of 
$\CP^{2}$. We may think of this as the region 
$$
\{(z_{1}, z_{2})\in \C^{2}: a^{2}\le |z_{1}|^{2} + |z_{2}|^{2}\le 1\}
$$
with boundaries collapsed along the Hopf flow, and give it the 
corresponding symplectic form $\om_{a}$.  
We are considering the action
$$
(z_{1}, z_{2})\mapsto (e^{-2\pi it}z_{1},e^{-2\pi it} z_{2}), \quad 0\le t \le 1.
$$
It is not hard to check that its normalized Hamiltonian  is
$$
H =  \pi (c-|z_{1}|^{2} - |z_{2}|^{2}), \qquad c = 
\frac{2(1 -a^{6})}{3(1 - a^{4})}.
$$
Since $ \max H = \pi(c -  a^{2})$ and $-\min H = \pi(1-c)$, 
we find $\max H > - \min H $ whenever $a^{2} < 1$.
(A similar example was given in Example 1.C in~\cite{Pac}.)  

The next task is to calculate 
$QH_{*}(M_{*})$.  This ring
is generated over $\La$ by elements $p = pt$, the exceptional 
divisor $E$, the fiber class $F = L-E$ and the fundamental class $[M_*]$. 
(Here $L = [\CP^1]$.)  
The quantum multiplication has $\1 = [M_*]$ as a unit, and is
derived from the following nontrivial Gromov--Witten invariants:
\begin{eqnarray*}
n(p, p, F; E+F) = 1;&& n(p, E, E; F) = 1;\\
n(A_{1}, A_{2}, A_{3}; E) = \pm 1&&\mbox{where }\;\; A_{i} = E\mbox{ 
or }F.
\end{eqnarray*}
One finds
$$
\begin{array}{lclcclcl}
 p*p &=& (E+F)\otimes e^{-E-F} &&& E * p  &=& F \otimes e^{-F}\\
p*F &=& \1\otimes e^{-E-F},\ &&&
E * E  &=& -p + E \otimes e^{-E} + \1\otimes e^{-F} \\
E * F &=&p  -E\otimes e^{-E} &&& F*F &=& E\otimes e^{-E}.
\end{array}
$$
We will be particularly interested in the element 
$Q = F\otimes e^{E/2+ F/4}$, since, as we shall see, this is the part 
of $\Psi(\al)$ that is independent of the choice of symplectic form 
$\om_a$.   Note that
$$
Q*Q = Q^2 = E\otimes e^{F/2},\qquad Q^{-1} = p\otimes e^{3F/4 + E/2},
$$
where the multiplication is $*$.  Recall that $v$ is defined by:
$$
v(\sum_{B\in \Hh} a_{B}\otimes e^{B}) = \sup \{\om(B): a_{B}\ne 0\}.
$$

\begin{lemma}\label{le:qk} 
 $v(Q^{k}) + v(Q^{-k}) \ge \om(F)$ for all $k>1$.
\end{lemma}
\proof{} 
First consider  $v(Q^{-k}), k\ge 1$.  The first few terms are
$$
\begin{array}{ccl}
Q^{-1}  & = & p\otimes e^{3F/4 + E/2},\\
Q^{-2} & = &  (E+F)\otimes e^{F/2},\\
Q^{-3} & = &   F \otimes e^{F/4 + E/2}  + \1 \otimes e^{F/4 - E/2},\\
Q^{-4} & = &  p \otimes e^{F} + \1.
\end{array}
$$
We claim that for all $m\ge 1$ and $1\le i \le 4$,
\begin{equation}\label{eq:coeff1}
    v(Q^{-i - 4m}) \ge   v(Q^{-i}).
\end{equation}
To see this note that
multiplication by $p$ has the effect:
$$
\begin{array}{cccc} p & E & F & \1\\
    \downarrow &  \downarrow &  \downarrow &  \downarrow \\
   E+F & F & \1 & p
    \end{array}   
$$
Since all coefficients are positive, no terms can cancel.
The $p$-term  in $Q^{-k}$ contributes to the $p$-term in $Q^{-(k+4)}$
with unchanged valuation via the cycle
$p\mapsto E\mapsto F\mapsto \1\mapsto p$. 
  Similarly, the $F$-term in 
$Q^{-k}$ contributes to the $F$-term in $Q^{-(k+4)}$
with unchanged valuation. Hence result. \MS

Next consider $Q^{k}, k>0$. The first few terms are
$$
\begin{array}{lcl}
Q  & = & F\otimes e^{E/2+ F/4},\\
Q^2  & = &  E\otimes e^{F/2},\\
Q^{3} & = & p \otimes e^{3F/4 + E/2} - E\otimes e^{3F/4 - E/2},\\
Q^{4} & = & - p \otimes e^{F} + E \otimes e^{F - E}+ \1,\\
Q^{5} & = & p \otimes e^{5F/4 - E/2} - E \otimes e^{5F/4 - 3E/2} + F \otimes 
e^{E/2 + F/4} - \1 \otimes e^{F/4 - E/2}.
\end{array}
$$
Thus the lemma holds for $1\le i\le 4$ by inspection.  We will  write 
$$
Q^{k} = p \otimes \la_{k,p} + E \otimes \la_{k,E} + F  \otimes 
\la_{k,F} + \1 \otimes \la_{k,\1},
$$
where $\la_{k,\cdot}\in \La$, and $\la_{k,\cdot} = 0$ for $k \le 0$,
and will prove 
that for all $m\ge1$
$$
\begin{array}{llll}
    v(\la_{i+4m,E}) &\ge &  v(\la_{i,E}),& i = 2,3\\
    v(\la_{1+4m,F}) & \ge &   v(\la_{1,F}), &\\
     v(\la_{4+4m,\1}) & \ge  & v(\la_{4,\1}).&
\end{array}
$$
In view of~(\ref{eq:coeff1}), this will prove the lemma.

Multiplying by $F$ has the following effect:
$$
\begin{array}{cccc} p & E & F & \1\\
    \downarrow &  \downarrow &  \downarrow &  \downarrow \\
    \1 & p, -E & E & F
    \end{array}   
    $$
The term in $E$ in $Q^{k}$ comes either from the term in $E$ in 
$Q^{k-1}$ or from the term in $F$ in $Q^{k-1}$.  The latter contribution
traces directly back to the term in $E$ in $Q^{k-4}$ via multiplication by 
the term $\1$ in $Q^4$.
Hence its valuation is $v(\la_{k-4,E})$, while the 
valuation of the former term is $v(\la_{k-1,E}) + \om(F/4- E/2)$.
Hence
$$
v(\la_{k,E}) = \max\{v(\la_{k-4,E}), v(\la_{k-1,E}) + \om(F/4- E/2)\}.
$$
Therefore $v(\la_{k,E})\ge v(\la_{k-4,E})$ provided that no terms cancel.
Cancellations could theoretically occur since the transition $E\mapsto E$ 
occurs with a $-$ sign; however we claim that they do not.  To see 
this, it suffices to look at the terms that contribute to
$\la_{2+4m,E}$.  Since $Q^{2}$ is a multiple of $E$
each term that involves $E$ in $Q^{2+4m}$ must involve some number of 
the period $4$ cycle $E\mapsto p\mapsto \1\mapsto F\mapsto E$ 
interspersed with some $E\to -E$ transitions.  But there must always 
be an even number of these transitions if we are to arrive back at $E$
after $4m$ multiplications.  Therefore  all the 
contributions to the $E$ term in $Q^{2+4m}$ occur with positive 
sign.
This shows that the coefficients $\la_{k,E}$ have the claimed 
behavior.  Similar arguments prove the statements about the other 
coefficients.
\QED

Our next aim is to calculate $\Psi(\al)$.   This depends on  
 $\om_{a}$, since $\om_{a}(E) = a^{2}, \om_{a}(F) = 1-a^{2}$.
 The case when $\om_{a}(F-2E) = 1-3a^{2} = 0$ is special, since in this 
 case $\om_{a}$ is a multiple of $c_{1}$, i.e. $(M_{*},\om_{a})$ is 
 monotone. We will see that when $3a^{2}\ne 1$ then
 $\Psi_{\al} = Q e^{\de(F-2E)}$ for an appropriate
 constant $\de$.   Since $c_{1}(F-2E) = 0$ the constant 
 $\de = \de(a)$  is determined by the requirement that 
 $u_{\al}(\si_{\al}) = 0$.  In fact, the exact value of $\de$ is 
 irrelevant for our main argument since we just need to estimate
 $\ell^{+}(k\al) + \ell^{+}(-k\al)$: see Corollary~\ref{cor:exam2} below.

It turns out to be easier first to 
calculate $\Psi(\la)$ where $\la = 2\al$.\footnote
{
A simpler way of doing these calculations is 
developed in~\cite{MT}.}

\begin{lemma} If $1-3a^{2}\ne 0$ then 
$$
    \Psi(\la) = 
  E\otimes e^{2\de(F-2E) + F/2},\quad\mbox{where }\;\de = 
 \frac{ (1-a^{2})^{2}}{12(1+a^{2})(1-3a^{2})}.
 $$
\end{lemma}
\proof{}  
Since the circle action $(z_{1},z_{2})\mapsto 
(e^{-2\pi it}z_{1},e^{-2\pi it} z_{2})$
preserves the fibers of the projection $M_*\to S^{2}$, the space
$P_{\la}$ fibers over $S^{2}\times S^{2}$ with fiber $S^{2}$.
Let us denote the generators of $H_{2}(S^{2}\times S^{2})$ by
$A = [S^{2}\times pt], \si = [pt\times S^{2}]$.  Here we are thinking 
that the original fibration $\pi: P_{\la}\to S^{2}$ is the composite
with projection $pr_{2}$ to the second factor, so that lifts of $\si$ to 
$P_{\la}$ correspond to sections of the original fibration $\pi$.
It is then not hard to see that the fibration over
$S^{2}\times S^{2}$ has the form 
$$
\pi_{1}\; : \; P(L\oplus \C) \to S^{2}\times S^{2}
$$
where $c_{1}(L) = A + \si$.  Further, this fibration has
an obvious complex structure $J$ as well as two natural
$J$-holomorphic sections, $Z_{-}$  with normal 
bundle $L^{-1}$ and $Z_{+}$ with normal 
bundle $L$.   Note that $M_{*} = \pi_{1}^{-1}(S^{2}\times pt)$. 
Moreover the class of the exceptional divisor  $E$ in $M_{*}$ is represented 
by the intersection $Z_{-}\cap M_{*}$.  To check this, 
recall that $P_\la$ is made from the region $R_H^+$ above the graph of $H$ by 
collapsing the boundary $\Ga_H$ to a single fiber.  Moreover,
the set of points $F_{\max}$ where the Hamiltonian $H$ takes its 
maximum is precisely the exceptional divisor, while the corresponding 
set $F_{\min}$ where $H$ is a minimum is a complex line.  It is not 
hard to see that in this realization the sections $Z_+$ and $ Z_{-}$ 
corrrespond to  the 
intersection of  $R_H^+$ with the slices $F_{\min}\times
[0,1]\times \R$ and $F_{\max}\times [0,1]\times \R$, respectively, 
where $Z_+$ corresponds to $F_{\min}$ because this slice has larger 
volume (with respect to  $\Om_0^2$.)

A dimension count shows that the 
only section classes $\si_{{\la}} + B$ that contribute to $\Psi(\la)$
are those  with vertical Chern class $c_{vert}\le 0$.  
Let $\si_-$ denote a lift to $Z_-$ of the sphere  $[pt\times S^2]$,
Then $c_{vert}(\si_{-}) = -1$, and
so the only classes that could contribute  $\Psi(\la)$ are
$\si_{-}$ and $\si_{-} + E$. Since $c_{vert}(\si_{-} + E) = 0$,
this class would contribute with coefficient
$$
 n_{P_{{\la}}}([M_*], [M_*], pt\,; \si_{-} + E).
 $$
 But this invariant is zero since all the $J$-holomorphic sections in
 class $\si_{-} + E$ lie in $Z_{-}= S^2\times S^2$ 
 and so do not meet an arbitrary 
 point.  (This follows by positivity of intersections since
 $Z_{-}$ is holomorphic and has \lq\lq negative'' normal bundle.)  
For similar reasons, the
 $J$-holomorphic sections in
 class $\si_{-}$ also lie in $Z_{-}$ and it is not hard to see that
 $$
 n_{P_{{\la}}}([M_*], [M_*], c\,; \si_{-}) = E\cdot c.
 $$

 It remains to calculate the canonical section class 
 $\si_{\la}$. We claim that
$$
\si_{\la} \; = \; \si_{-} + F/2 + 2\de(F-2E)
$$
for the  $\de$ given above.
First observe that, as in~\cite{Mcq}~Example~3.1, the cohomology ring
$H^{*}(P_{\la})$ is generated by the vertical Chern class $\nu$ of 
$\pi_{1}: P\to S^{2}\times S^{2}$ together with the pullbacks 
$\mu_{1}, \mu_{2}$
via $\pi_{1}$ of the two obvious generators of $H^{2}(S^{2}\times 
S^{2})$.  Thus
$$
\mu_{1}(E) = 1, \;\;\nu(E) = -1,\;\; \mu_{1}(F) = \mu_{2}(F) = 0,\;\;
\nu(F) = 2.
$$
Clearly $\mu_{i}^{2}= 0$.  It is also not hard to see that 
 $\nu^{2} = 2\mu_{1}\mu_{2}$. (As in~\cite{Mcq}~Example~3.1, look at 
the  Poincar\'e duals.) 
The vertical Chern class $c_{vert}$ for $P_{\la}\to S^{2}$  takes the 
value $1$ on $E$ and so is
$\nu + 2\mu_{1}$. But $\mu_{1}(\si_{-}) = 0$. Therefore
$c_{vert}(\si_{-} + F/2) = 0$.  Since $c_{vert}(F - 2E) = 0$ and
$u_{\la}(F - 2E) = \om_{a}(F - 2E) \ne 0$ when $3a^{2}\ne 1$,
$\si_{\la}$ must have the stated form for some $\de$.
Now observe that 
$$
u_{\la} = \frac{1-a^{2}}2 \nu + \frac {1+a^{2}}2 \mu_{1} + \eps \mu_{2},
$$
and use the fact that $u_{\la}^{3} = 0$ to conclude $\eps = 
-(1-a^{2})^{2}/6(1+a^{2})$.  Further, because $u_{\la}(\si_{\la}) = 0$
we must have $\eps = 2\de(3a^{2}-1)$, as required.

Finally, the fact that  $\si_{-} = \si_{\la} - F/2 - 2\de(F-2E)$ 
implies
$$
\Psi({\la}) = E\otimes e^{F/2 + 2\de(F-2E)}
$$
as claimed.
(Compare Seidel~\cite{Seid}, where a  similar calculation 
is made 
for $M_* = S^{2}\times S^{2}$.)\QED

\begin{prop}\label{prop:q} If $3a^{2}\ne 1$, then $\Psi({\al}) = 
 \pm F\otimes e^{E/2  + F/4 + \de(F-2E)}$ where $\de$ is as before.
\end{prop}
\proof{} 
We show that
$\Psi(\al)$ has the form $E\otimes \la_E + F\otimes \la_F$,
i.e. that the coefficients of $p$ and $\1$ vanish.  The desired result then 
follows by an easy calculation, using our knowledge of
$(\Psi(\al))^2$.

The action
$$
(z_1,z_2)\mapsto (e^{-2\pi i\theta}z_1,z_2)
$$
permutes the lines through $(0,0)\in \C^2$. 
 Hence its action 
on $M_*$ is the lift of  an action $\al^*$ on the base $S^2$ of the obvious 
fibration  $M_* \to S^2$.  It is easy to see that $\al^*$ is a 
rotation with two fixed points, one corresponding to the line $z_1 = 0$ 
(whose points are all fixed) and the other corresponding to 
$z_2 = 0$.  Hence the fibration $\pi: P_\al \to S^2$ factors as
$$
P_\al \stackrel{\pi_1}\to P_{\al*} \stackrel{\pi_2}\to S^2,
$$
where $P_{\al^*}$ is a one point blow up of $S^2$ and 
$\pi_1:P_\al \to P_{\al^*}$
is a fibration with fiber $S^2$.   All these spaces have natural 
complex structures that are preserved by these maps.   Further there are 
two natural sections $Z_\pm$  of $\pi_1$ given by the restricting the 
action $\al$ to the two canonical sections of $M_*\to S^2$ (the images of 
the two boundary components of the region $\{a^2 \le |z_1|^2 + |z_2|^2 \le 
1\}$.) 

The fibration
$\pi_2:P_{\al^*} \to S^2$ has two natural sections given by  the two
fixed lines in $M_*$.  
The exceptional divisor $E^*$  in $P_{\al*}$  corresponds to 
the line $z_1 = 0$ which is pointwise fixed by $\al$, while the other 
section $L^*$ has self intersection $+1$.  
It follows that the fibration
$$
\pi: \pi_1^{-1}(E_*) \to S^2
$$
is holomorphically trivial: it has the form $S^2\to P_{\ga}\to S^2$ where
$\ga$ is the induced action on the line $z_1 = 0$. On the other hand the 
fibration  $\pi: \pi_1^{-1}(L_*) \to S^2$ is nontrivial, and its total 
space is  the one point blow up of $\CP^2$.  Indeed given any 
other line $L\in P_{\al^*}$, the total space $\pi_1^{-1}(L)$ is also the 
one point blow up of $\CP^2$  with exceptional divisor 
$Z_-\cap (\pi_1)^{-1}(L)$.  

Now observe that the coefficient of $p$ in $\Psi(\al)$ is nonzero 
only if there is a nonzero invariant of the form
$n_{P_\al}([M_*], [M_*], [M_*]; \si)$.  Hence $c_{vert}(\si) = -2$
and $\si$ is represented (for generic $J$ on $P_\al$) by isolated curves.
The natural complex structure on $J_\al$ may not be regular for $\si$.
Nevertheless, since it is a limit of generic structures, Gromov 
compactness implies that it would have to 
contain some holomorphic representative of $\si$ which would be a section 
$S$
together perhaps with some bubbles in the fibers.  Since each fiberwise 
bubble has positive Chern class, this means that  
$c_{vert}(S) \le -2$.  But  such a section $S$ does 
not exist.  It would have to project to a section $S_2$ of $\pi_2:P_{\al^*}\to 
S^2$.  We cannot have $S_2 = E_*$ since the only lifts of $E_*$ to $P_\al$ 
have $c_{vert} = -1$. Thus $S_2$ would have to be a line $L$ in 
$P_{\al^*}$,
and, because $c_{vert}(L) = 1$,
 $S$ would have to have self intersection $-3$ in $\pi_1^{-1}(S_2)$.  
This is impossible since, as we saw above, $\pi_1^{-1}(S_2)$ is the blow up 
of $\CP^2$.

A similar argument shows that the coefficient of $\1$ in $\Psi(\al)$ must 
vanish.  Otherwise there would be a section class $\si$ in $P_\al$
with $c_{vert} = 0$ 
with holomorphic representatives through every point $x\in P_\al$.  Again these 
representatives would have to be the union of a section $S'$ together with 
some fiberwise bubbles.  Hence if $\pi_1(x)\notin E_*$, $S'$ would have 
to project to a line $L$  in $P_{\al^*}$.  Such a line has one lift to a 
curve with $c_{vert} = 0$, namely $\pi_1^{-1}(L)\cap Z_-$, but 
 such a lift does not exist through an arbitrary point.\QED

\begin{cor}\label{cor:exam2}  $r_{\trho}(M_{*}, \om_{a}) = 
    r_{\trho_{f}}(M_{*}, \om_{a}) = \om(F) = 
\pi(1-a^{2})$.
\end{cor}
\proof{} $\pi_{1}(\Ham(M_{*}, \om))$ is isomorphic to $\Z$ with 
generator $\al$.   It follows from Corollary~\ref{cor:main}
and Lemma~\ref{le:qk} that  when $3a^{2}\ne 1$
\begin{eqnarray*}
\trho(k\al) \;\ge\; \trho_f(k\al)& = &\trho^+(k\al) + \trho^+(-k\al) \\
& \ge & v(Q^k e^{k\de(F-2E)}) + v(Q^{-k}e^{-k\de(F-2E)})\\
& = & 
v(Q^k) + v(Q^{-k})  \;\ge\; \om(F),
\end{eqnarray*}
where the second equality holds because the terms involving $\de$ 
cancel.
But the loop $2\al$ is generated by the Hamiltonian
$$
H = \pi(c - |z_1^{2}| - |z_2^{2}|) 
$$
which has length $\pi(1 - a^{2}) = \om(F).$   Hence result.

This conclusion is still valid when $3a^{2} = 1$ since the
 inequality $\trho(k\al) \ge \om(F)$ still holds.  Because it is uniform 
 one can prove this by continuity in $a$: the set of $a$ for which 
 this inequality does not hold must be open.  Alternatively, one can argue
 using the remarks 
 after Corollary~\ref{cor:main}.  
\QED

\NI
{\bf Proof of Proposition~\ref{prop:exam2}}.\,\,
The above corollary shows that the loop $2\al$ on $(M_{*}, \om_{a})$ satisfies 
all the required conditions.\QED

\begin{rmk}\rm\begin{itemize}\item[(i)] It follows from~\cite{MSlim,En} that the 
circle action generated by $K = -\pi|z_1^{2}|$  is also a 
$\trho$-minimizing representative of its homotopy class $\al$.
Thus 
$$
\ell_{\trho} (\al) = \pi \;\;>\;\; \ell_\trho(2\al) = \pi(1-a^{2}).
$$ 
\item[(ii)]  It is also interesting to try to understand the
asymptotic growth of the loops $\pm \al$ by calculating the limits
of the one sided measures $\ell^+ (k\al)/k$ 
as $k\to \pm \infty$ for all values of $a< 1$.
Polterovich showed in~\cite{Pac} that in the monotone case 
$\ell^+(k\al)/k\ge |I(\trho)| = 1/18$ when $k\ge 1$ and $I$ is as in 
Remark~\ref{rmk:pol}, but his method gave nothing for $k\to -\infty$.
Our methods extend his result for the case $k\to \infty$ to all values of 
$a$, but also have no information about the limit as  $k\to -\infty$.
To see this, note first that
$\om(\de(F-2E)) > 0$ for all $a$. 
The arguments in Lemma~\ref{le:qk} show that
$v(Q^{-k})$ is bounded as $k\to \infty$ when $3a^{2}\ge 1$ and grows as
$k\om(F/4 - E/2)/3$ for $3a^{2} < 1$.  (This growth comes from the 
cycle $p\mapsto F\mapsto \1\mapsto p$ that increases 
valuation by $\om(F/4 - E/2)$.)  Since $\om(F/4 - E/2)/3 < \om(\de(F-2E)) $
for all $a$, $v(\Psi(-k\al)) \to 
-\infty$ as $k\to \infty$.  Thus we get no information on 
$\ell^+(-k\al)$.  On the other hand, $v(Q^{k})$ is either 
bounded, or, if $3a^{2} < 1$, grows as a multiple of $\om(F/4-E/2)$.
Hence 
$$
\ell^+(k\al)/k \;\ge\; v(\Psi(k\al))/k\; \ge \om(\de(F-2E)) = 
 \frac{ (1-a^{2})^{2}}{12(1+a^{2})}\; \quad\mbox{as }\;\; 
k\to \infty
$$
for all $a$.
 See Polterovich~\cite{Pbk}
for a further discussion of this question.
\end{itemize}
\end{rmk}

  \end{document}